\documentclass[11pt]{amsproc}

\usepackage{hyperref} 

\usepackage[all]{xy} 

\usepackage{txfonts}

\newtheorem{theorem}{Theorem}[section]
\newtheorem{lemma}[theorem]{Lemma}
\newtheorem{proposition}[theorem]{Proposition}
\newtheorem{corollary}[theorem]{Corollary}

\theoremstyle{definition}
\newtheorem{definition}[theorem]{Definition}

\theoremstyle{remark}
\newtheorem{remark}[theorem]{Remark}

\numberwithin{equation}{section}

\newcommand{\abs}[1]{\lvert#1\rvert}

\newcommand{\C}{\ensuremath{\mathbb{C}}}

\newcommand{\Q}{\ensuremath{\mathbb{Q}}}
\newcommand{\R}{\ensuremath{\mathbb{R}}}
\newcommand{\Z}{\ensuremath{\mathbb{Z}}}

\newcommand{\eps}{\ensuremath{\varepsilon}}
\newcommand{\HH}{\ensuremath{\mathbb{H}}}

\newcommand\al{\alpha}
\newcommand\bt{\beta}
\newcommand\g{\gamma}
\newcommand\dt{\delta}
\newcommand\e{\varepsilon}
\renewcommand\th{\vartheta}

\newcommand\x{\xi}
\newcommand\ch{\chi}
\newcommand\ps{\psi}
\newcommand\Om{\Omega}
\newcommand\om{\omega} %\newcommand{\HH}{\ensuremath{\mathfrak{H}}}

\newcommand{\cK}{\ensuremath{\mathcal{K}}}
\newcommand{\cV}{\ensuremath{\mathcal{V}}}

\newcommand{\cH}{\ensuremath{\mathrm{Fib}}}
\newcommand{\cL}{\ensuremath{\mathcal{L}}}
\newcommand\gr{\phi}

\newcommand{\FEd}{\ensuremath{{}^3\mathrm{FE}}}
\newcommand{\FEv}{\ensuremath{{}^4\mathrm{FE}}}

\newcommand{\TT}{\ensuremath{T^\prime}}
\newcommand{\dd}{\ensuremath{\mathrm{d}}}
\newcommand{\Gmod}{\ensuremath{\Gamma}}

\newcommand\fs{{\om^\ast}}
\newcommand\smp{{\om^\ast\!,\mathrm{simple}}}
\newcommand\sa{{\om^\ast\!,\infty}}
\newcommand\parb{{\mathrm{par}}}

\newcommand{\re}[1]{\ensuremath{{\mathrm{Re}(#1)}}}
\newcommand{\im}[1]{\ensuremath{{\mathrm{Im}(#1)}}}

\newcommand{\GL}[1]{\ensuremath{{\mathrm{GL}\!\left(2, #1 \right)}}}
\newcommand{\PSL}[1]{\ensuremath{{\mathrm{PSL}\!\left(2, #1 \right)}}}
\newcommand{\PGL}[1]{\ensuremath{{\mathrm{PGL}\!\left(2, #1 \right)}}}
\newcommand{\rMatrix}[4]{{\begin{pmatrix} #1 & #2 \\ #3 & #4
\end{pmatrix}}}
\newcommand{\Matrix}[4]{{\begin{bmatrix} #1 & #2 \\ #3 & #4
\end{bmatrix}}}
\newcommand{\LL}{\ensuremath{\cL}}
\newcommand{\LLt}[1]{\ensuremath{{\tilde\cL}_{#1}}}
\newcommand{\LLMayer}[1]{\ensuremath{\cL^\mathrm{Ma}_{#1}}}
\newcommand{\sing}[1]{\ensuremath{\mathrm{Sing}\left({#1}\right)}}
\newcommand{\PPP}[1]{\ensuremath{\mathbb{P}^1_{#1}}}
\newcommand{\aver}{\ensuremath{\av_T}} % average 
\newcommand{\avp}{\ensuremath{\aver^+}} % right side average of T
\newcommand{\avm}{\ensuremath{\aver^-}} % left side average of T
\newcommand{\avpm}{\ensuremath{\aver^\pm}}
\newcommand\av{\mathrm{A\!v}}
\newcommand{\OO}[1]{\ensuremath{\mathcal{O}\left( #1 \right)}}

\newcommand{\PP}{\PPP{\R}}
\newcommand\txtfrac[2]{{\textstyle \frac{#1}{#2}}}
\renewcommand\setminus\smallsetminus

\begin{document}

\title[Transfer operators and cohomology]{Eigenfunctions of transfer
operators and cohomology}
\date{July 4, 2007}

%    Information for first named author
\author{R.W.~Bruggeman}
\address{\href{http://www.math.uu.nl}{Mathematisch Instituut},
\href{http://www.uu.nl}{Universiteit Utrecht}, Utrecht, The
Netherlands}
\email{\href{mailto:bruggeman@math.uu.nl}{bruggeman@math.uu.nl}}

%    Information for second named author
\author{T.~M\"uhlenbruch}
\address{\href{http://www.dynamik.tu-clausthal.de/}{Institut f\"ur
Theoretische Physik}, \href{http://www.tu-clausthal.de}{Technische
Universit\"at Clausthal}, Clausthal-Zeller\-feld, Germany}
\curraddr{Mathematical Biology Unit, Okinawa Institute of Science and
                Technology, 7542 Onna, Onna Village, Okinawa 904-0411, Japan}
%\email{\href{mailto:tobias.muehlenbruch@tu-clausthal.de}%
%{tobias.muehlenbruch@tu-clausthal.de}}
\email{\href{mailto:muehlenbruch@oist.jp}{muehlenbruch@oist.jp}}
\thanks{The second named author was supported in by the
\href{http://www.dfg.de/}{Deutsche Forschungsgemeinschaft} through
the DFG Research Project ``Transfer operators and non arithmetic
quantum chaos'' (Ma 633/16-1).}

%    General info



\begin{abstract}%text is identical to that in file tro_abstract
% except for $-signs around numbers
The eigenfunctions with eigenvalues $1$ or $-1$ of the transfer
operator of Mayer are in bijective correspondence with the
eigenfunctions with eigenvalue $1$ of a transfer operator connected
to the nearest integer continued fraction algorithm. This is shown by
relating these eigenspaces of these operators to cohomology groups
for the modular group with coefficients in certain principal series
representations.
\end{abstract}

\maketitle
%%%%%%%%%%%%%%%%%%%%%%%%%%%%%%%%%%%%%

%%%%%%%%%%%% % Table of contents

\setcounter{tocdepth}{1} \tableofcontents

%%%%%%%%%%%% Section 1 %%%%%%%%%%%%%%%

\section{Introduction} \label{A}

D.Mayer, \cite{Ma90}, defined the operator
\begin{equation}
\label{A.1}
\LLMayer{s}f(z) = \sum_{n=1}^\infty (z+n)^{-2s} \, f\left(
\frac{1}{z+n} \right)
\end{equation}
on the Banach space of continuous functions on the disk
$\abs{z-1} \leq \frac{3}{2}$, holomorphic on
$\abs{z-1} < \frac{3}{2}$, with the supremum norm. The series
converges absolutely if $\re{s} > \frac{1}{2}$. There is a
meromorphic continuation in $s$, with a pole at $\frac12$ as the sole
singularity in the region $\re s>0$. The operator $\LLMayer s$ is a
transfer operator for the Artin billiard dynamical system \cite{Ar}.
It is connected to the Gauss map $x\mapsto \frac1x-[\frac1x]$ on
$[0,1]$. Ultimately, this dynamical system comes from closed billiard
flows on the quotient of the upper half plane by $\PGL\Z$. (Here
$\pm \rMatrix{-1}001$ acts by $z\mapsto -\bar z$.)
An introductory lecture on Mayer's transfer operator is~\cite{Ma}.

The DFG research project \emph{Transfer operators and non arithmetic
quantum chaos} (Ma 633/16-1) involves finding transfer operators
connected to the dynamical systems of closed geodesic flows on the
hyperbolic surfaces represented by the quotient of the upper half
plane by arbitrary Hecke triangle groups. For the modular group, it
leads to another transfer operator $\LLt{s}$. This operator acts in
the space of vectors of two holomorphic functions on the open unit
disk which are continuous on the closed unit disk. With the supremum
norm these vectors form a Banach space. The operator is given by
\begin{eqnarray}
\label{A.5}
\LLt s \vec f
&=& \rMatrix{\tilde \LL_{s}^{1,1}}{\tilde \LL_{s}^{1,2}}{\tilde
\LL_{s}^{2,1}}{\tilde \LL_{s}^{2,2} } \vec{f} \qquad\mbox{with }\\
\nonumber
\LLt s ^{1,1} f_1(z)
&=& \sum_{n=3}^\infty (z+n)^{-2s} f_1\left( \frac{-1}{z+n}\right), \\
\nonumber
\LLt s ^{1,2} f_2(z)
&=& \sum_{n=2}^\infty (n-z)^{-2s} f_2\left( \frac{1}{-z+n}\right), \\
\nonumber
\LLt s ^{2,1} f_1(z)
&=& \sum_{n=2}^\infty (z+n)^{-2s} f_1\left( \frac{-1}{z+n}\right),
\qquad \mbox{and}\\
\nonumber
\LLt s ^{2,2} f_2(z)
&=& \sum_{n=3}^\infty (n-z)^{-2s} f_2\left( \frac{1}{-z+n}\right)
\end{eqnarray}
This converges absolutely for $\re s>\frac{1}{2}$, and has a
meromorphic continuation in~$s$ with a first order pole at
$s=\frac12$ as sole singularity in the region $\re s>0$.

Our main result is:
\begin{theorem}\label{thm-main}Let $s\in \C$, $0<\re s<1$,
$s\neq\frac12$. There is a bijective correspondence between the
spaces
$\ker(\LLMayer s-\nobreak1)  \oplus \ker(\LLMayer s+\nobreak1)$ and
$\ker(\LLt s - \nobreak1)$.
\end{theorem}

The eigenfunctions of both transfer operators satisfy finite linear
identities. Lewis and Zagier \cite[Proposition in \S3,
Chap.~IV]{LZ01} show that if $\LLMayer s f=\pm f$, then
$P(z)=f(z-\nobreak1)$ satisfies the three term equation
\begin{equation}
\label{3te}
P(z) \;=\; P(z+1)+(z+1)^{-2s}P\left( \frac z{z+1}\right)
\end{equation}
and the parity condition
\begin{equation}
z^{-2s}P(1/z) \;=\; \pm P(z)\,.
\end{equation}
These functions extend holomorphically to
$\C'=\C\setminus(-\infty,0]$. So $\ker(\LLMayer s-\nobreak1) 
\oplus \ker(\LLMayer s+\nobreak1)$ corresponds to a subspace of the
space of all holomorphic solutions of \eqref{3te} on~$\C'$. This
subspace is characterized by the asymptotic behavior
$P(x) = c_\infty x^{1-2s}+O(x^{-2\re s})$ as $x\uparrow\infty$, and
$P(x) = c_0 x^{-1} + O(1)$ as $x\downarrow0$.

For both transfer operators, we relate, in \S\ref{C}, the
eigenfunctions on disks to eigenfunctions in the real analytic
functions on an interval. This allows a cohomological interpretation,
to be discussed in \S\ref{sect-efMa} and \S\ref{D}. We will show that
solutions of $\LLt s \vec f=\vec f $ correspond to solutions of the
four term equation
\begin{equation}\label{4te}
g(z) + (z+2)^{-2s} g\left(\frac{-1}{z+2}\right) \;=\; g(z-1) +
(2-z)^{-2s} g\left( \frac{1-z}{z-2}\right)\,,
\end{equation}
on a suitable domain containing $(-1,1)$.

Not all solutions of \eqref{3te}, respectively \eqref{4te}, correspond
to eigenfunctions of $\LLMayer s$, respectively eigenfunctions of
$\LLt s$. We will see in Theorem~\ref{thm-th} and
Proposition~\ref{prop-an-coh} that the space of all real analytic
solutions of \eqref{4te} on the interval
\[\left( \frac{-3-\sqrt5}2,\frac{1+\sqrt 5}2\right)\]
is isomorphic to the first cohomology group of the modular group with
coefficients in the principal series representation with spectral
parameter~$s$. Theorem~\ref{thm-par} shows that the eigenspace of
$\LLt s$ for the eigenvalue~$1$ corresponds to a well defined
subspace of this cohomology group. This same subspace is also
isomorphic to the sum of the eigenspaces of $\LLMayer s$ for the
eigenvalues $1$ and $-1$. This can be shown by methods in \cite{LZ01}
and~\cite{BLZ}.

We will take care to indicate the various steps in the correspondence
between eigenspaces of $\LLMayer s$ and $\LLt s$ as explicitly as
possible, even for steps where we might refer to \cite{LZ01}
or~\cite{BLZ}. The least explicit step is an application of
Proposition~\ref{D2.6}, where a function with two singularities is
written as a difference of two functions which each have a
singularity in only one point.
\medskip

As background information, we discuss in \S\ref{B8} how $\LLt s$
arises from the nearest integer continued fraction transformation.
Our proof of the correspondence does not use that both transfer
operators arise from the geodesic flow on related quotients of the
upper half plane. It would be interesting to go directly from the
geodesic flow to the relevant cohomology groups.

The Selberg trace formula relates the recurrent points of the geodesic
flow to spectral data. So both transfer operators have a relation to
Maass forms.

Our cohomological approach to the correspondence is based on ideas in
\cite{LZ01} and \cite{BLZ}. The leading idea in \cite{BLZ} is the
relation between certain cohomology groups and spaces of Maass forms,
which we discuss in \S\ref{B7}. This relates eigenfunctions of
$\LLMayer s$ and $\LLt s$ to Maass forms without use of the Selberg
trace formula.\medskip

We thank R.~Sinclair for his remarks on a preliminary version of this
paper. We thank F.~Str\"omberg for the fruitful discussions of the
        transfer operator and the underlying dynamical system.

\subsection{The transfer operator for the nearest integer continuous
fraction algorithm} \label{B8}
Although the origin of $\LLt s$ from a dynamical system is not used in
this paper, it seems right to explain why $\LLt{s}$ deserves to be
called a \emph{transfer operator}. It is derived by the Ruelle
 transfer operator method applied to a dynamical system based on the
nearest integer map and associated continued fractions. Such nearest
integer continued fractions have already been discussed in~\cite{Hu}.

Consider the interval map
\begin{equation}
\label{B8.3}
f_3: \quad I_3 \to I_3 ; \qquad x \mapsto S\,x - \left\lfloor S\,x +
\frac{1}{2}\right\rfloor
\end{equation}
with $I_3 =\left[ -\frac{1}{2}, \frac{1}{2} \right]$ and
$\lfloor x \rfloor$ the element $n\in \Z$ such that $n\leq x <n+1$ if
$x>0$, and $n<x\leq n+1$ if $x\leq0$. The function $f_3$ is closely
related to the nearest integer continuous fractions.
(We keep the the subscript $3$ since it is a specialization of an
interval map associated to Hecke triangle groups.)
Basically $f_3$ acts as the ``left-shift'' on the space of
configurations $(a_1,a_2,\ldots)$ for the nearest integer continued
fraction expansion
\[ [0;a_1,a_2,\ldots]:= \frac{-1}{a_1 + \frac{-1}{a_2 +
\frac{-1}{\ldots}}} \in I_3\,. \]
The map $f_3$ generates a discrete dynamical system of finite type.
The transfer operator associated to $f_3$ is
\[ T_s \, f(x) = \sum_{y \in f^{-1}_3 (x)} \left|\frac{\dd
f_3^{-1}(x)}{\dd x} \right|^s \; f(y)
\]
defined on a suitable function space. The expression
$\frac{\dd f_3^{-1}(x)}{\dd x}$ denotes the derivative of the
appropriate invertible branch of $f_3^{-1}$ in $x$. We find
$\LLt{s} = T_s$ on the Banach space on which $\LLt{s}$ is defined.
Forthcoming work in the project Ma633/16-1 of the Deutsche
Forschungsgemeinschaft will give more details.
\smallskip

The dynamical system is related to the geodesic flow on the hyperbolic
surface $\PSL \Z\backslash \HH$. One can show that the
Fredholm-determinant $\det \left( 1 - \LLt{s} \right)$ is essentially
equal to the Selberg zeta-function $Z(s)$ for the full modular group
$\Gmod$. Here, essentially equal means that
\begin{equation}
\label{B8.1}
Z(s) = \frac{\det \left( 1 - \LLt{s} \right)}{\det \left( 1 - \cK_s
\right)}
\end{equation}
with $\cK_s$ defined on the same space of pairs of function by
\begin{equation}
\label{B8.2}%B6.1
\cK_s { g_1 \choose g_2} = {g_1 \big|_{2s} ST^3 \choose g_1 \big|_{2s}
ST^3}.
\end{equation}
It is known that $\cK_s$ only admits the eigenvalue~$1$ if
$s \in \{ -n + \frac{\pi ik}{2\ln (r)}; \; 
n \in \Z_{\geq 0}, k \in \Z\}$. These values of~$s$ are not in the
domain under consideration in this note.

The zeros of $Z(s)$ on the line $\re{s} = \frac{1}{2}$ correspond to
eigenvalues $s(1-\nobreak s)$ of the hyperbolic Laplace operator.
  Hence it is not surprising that eigenfunctions of~$\LLt s$ with
eigenvalue~$1$ can be related to cohomology classes that are
themselves related to Maass forms. See~\S\ref{B7}.
 \medskip

The advantages of $\LLt{s}$ compared to Mayer's transfer operator
$\LLMayer{s}$ are that its construction is directly related to the
geodesic flow and that the same construction works for all Hecke
triangle groups. The disadvantage of $\LLt{s}$ compared to
$\LLMayer{s}$ is its more complicated structure.

%%%%%%%%%%%% Section 2 %%%%%%%%%%%%%%%

\section{Definitions and results} \label{B}
In this preliminary section we define or recall various concepts to be
used in this paper. Among them are the principal series
representations in \S\ref{B3}, and the definition of parabolic
cohomology groups in \S\ref{B4}.  We state in \S\ref{sect-LM-coh} and
\S\ref{sect-hypcoh} the results implying the main result
Theorem~\ref{thm-main}. Finally we give in \S\ref{B7} background
information on the relation between certain spaces of Maass forms and
cohomology groups.

\subsection{Modular group} \label{B1}
We use $\Matrix abcd$ to denote
$\left\{ \rMatrix {ta}{tb}{tc}{td} \;:\; t\neq0\right\}$ in
\[ \PGL\R = \GL \R / \left\{ \rMatrix t00t \;:\; t\in
\R\setminus\{0\}\right\}\,.\]

We work with the full modular group $\Gmod=\PSL \Z$, which is discrete
in $\PSL\R$, and is generated by $S=\Matrix 0{-1}10$ and
$T=\Matrix 1101$, with the relations
\begin{equation}
\label{B1.1}
S^2 = (ST)^3 = 1.
\end{equation}
We denote also
\begin{equation}
\label{B1.3}
\TT : = TST = ST^{-1}S= \Matrix{1}{0}{1}{1} \in \Gmod
\quad \mbox{and} \quad C = \Matrix{0}{1}{1}{0} \in \PGL{\Z}.
\end{equation}
We have $\PGL{\Z} = \Gmod \cup C \, \Gmod$.

\subsection{Principal series representations} \label{B3} We
describe the standard
realization of the principal series representations of $\PGL\R$ in
the functions on~$\R$.

The group $\PGL\R$ has a family of actions, parametrized by $s\in \C$,
on functions defined on a subset of~$\R$, given by
\begin{equation}\label{B2.1a} h|_{2s} M \,(x) = \abs{ad-bc}^s \,
\abs{cx+d}^{-2s} \, h\!\left(\frac{ax+b}{cx+d}\right). \end{equation}
This is a right action; the natural place for the symbol $|_{2s}$ is
after the function. We call $s$ the \emph{spectral parameter}.

For each value of~$s$, this action preserves the spaces $\cV^\omega_s$
and $\cV^\infty_s$ of real-analytic and smooth vectors in the
discrete series representation with spectral parameter~$s$. The space
$\cV^\omega_s$ consists of the $h:\R\rightarrow\C$ that are
real-analytic on~$\R$ and for which $x\mapsto |x|^{-2s}h(-1/x)$
extends as a real-analytic function on a neighborhood of~$0$. Replace
`real-analytic' by `smooth' to obtain the characterization of
$h\in \cV^\infty_s$.

We refer to \cite[Section~2]{BLZ} for a discussion of other models of
the principal series representations. Here it suffices to note that
elements of $\cV^\omega_s$ and $\cV^\infty_s$ can be viewed as
functions on the projective line $\PP=\R\cup\{\infty\}$, and that the
required behavior of
\begin{equation}
\label{B3.1}
x\mapsto \left(h|_{2s} S\right)(x) =|x|^{-2s} h(-1/x)
\end{equation}
may be viewed as the description of analyticity or smoothness
at~$\infty$.

The real-analytic functions in $\cV^\om_s$ are the restriction of a
holomorphic function on some neighborhood of~$\R$, that depends on
the functions. On such holomorphic functions the slash operator takes
the form
\begin{equation}
\label{B2.1}
h|_{2s} M \,(z) := \abs{ad-bc}^s \, \big((cz+d)^2\big)^{-s} \,
h\!\left(\frac{az+b}{cz+d}\right).
\end{equation}
The factor $\big((cz+d)^2\big)^{-s}$ is holomorphic in~$z$ for
$\re z\neq -\frac dc$. In some cases, for
instance for period functions satisfying \eqref{3te}, we may prefer
to choose the factor differently, such that it is holomorphic on the
domain of the function and positive on the real points in the domain.

We need more spaces related to $\cV^\omega_s$ and $\cV^\infty_s$. If
$I\subset\PP$ is an open subset, we define $\cV^\omega_s(I)$ as the
space of $h:I\cap\R\rightarrow \C$ that are real-analytic
on~$I\cap \R$, and for which in the case $\infty\in I$ the function
in \eqref{B3.1} is real-analytic at~$0$. In particular, if
$I=\PP\setminus E$ for some finite set $E$, then $\cV^\omega_s(I)$
consists of analytic vectors with finitely many singularities
on~$\PP$.

The space $\cV^{\fs}_s$ is defined as the inductive limit
\begin{equation*}
\cV^{\fs}_s = \lim_{\longrightarrow} \cV^\omega_s(\PP\setminus E),
\end{equation*}
where $E$ runs through the finite subsets of~$\PP$. If
$h\in \cV^{\fs}_s$, then there is a minimal finite set $E\subset\PP$
such that $h\in \cV^\om_s(\PP\setminus E)$. We call the elements of
$E$ the \emph{singularities} of~$h$ and denote this set by $\sing h$.

By imposing conditions at the singularities, we define subspaces of
$\cV^\fs_s$. For instance
\begin{equation}
\cV^{\sa}_s := \cV^{\fs}_s \cap \cV^\infty_s
\end{equation}
is the space of smooth vectors that are real-analytic outside finitely
many points. A slightly larger space is $\cV^{\smp}_s$ consisting of
the $h\in \cV^{\fs}_s$ for which we allow simple pole at finitely
 many points. Its elements $f$ have to satisfy at their finitely many
singularities $x_0$:
\begin{eqnarray}
\label{simple-sing}
&& x \mapsto (x-x_0) \, f(x) \mbox{ is smooth at } x_0 \mbox{ if } x_0
\not= \infty, \mbox{ and}\\
\nonumber
&& y \mapsto y \, \abs{y}^{-2s} \,f (-1/y) = y \, \left(f|S \right)(y)
\mbox{ is smooth at } 0 \mbox{ if } x_0 = \infty\,.
\end{eqnarray}

Thus we have various spaces, all invariant under the action $|_{2s}$
of $\PGL\R$, that satisfy the following inclusions:
\[ \xymatrix{ \cV^\om_s\ar@{^{(}->}@<-.4ex>[r]
& \cV^{\sa}_s\ar@{^{(}->}@<-.4ex>[r] \ar@{^{(}->}[d]
& \cV^{\smp}_s \ar@{^{(}->}@<-.4ex>[r]
& \cV^{\fs}_s\\
& \cV^\infty_s & } \]

Throughout this note we use the assumption
\begin{equation}
\label{B2.2}
0<\re s<1,\qquad s \neq \frac{1}{2}.
\end{equation}
Mostly, we work with a fixed value of the spectral parameter~$s$. Then
we shall write $h|M$ instead of $h|_{2s}M$.

For vector valued functions we write ${f_1\choose f_2} |_{2s} M$ for
${f_1 |_{2s} M \choose f_2 |_{2s} M}$. The slash operator extends to
elements of the group ring $\C[\PGL\R]$ by
\[ f|_{2s} (M_1 + a M_2) = f|_{2s} M_1 + a \, f|_{2s} M_2
\qquad\mbox{for all } g_1,g_2 \in \Gmod \mbox{ and } a \in \C. \]
In some circumstances one can make sense of $f|_{2s}\Xi$ where $\Xi$
is an infinite sum of elements of $\Gmod$.

\subsection{Cohomology groups} \label{B4}
As usual, the \emph{first cohomology group} of $\Gmod$ with values in
a right $\Gmod$-module $V$ can be described by
\begin{eqnarray}
\label{B4.1}
H^1(\Gmod; V)
&=& Z^1(\Gmod; V) / B^1(\Gmod; V),\\
\nonumber
Z^1(\Gmod; V)
&=& \big\{ \psi: \Gmod \to V; \; \psi_{\gamma \delta} = \psi_\gamma |
\delta + \psi_\delta \mbox{ for all } \gamma,\delta \in \Gmod \big\}
\quad \mbox{and}\\
\nonumber
B^1(\Gmod; V)
&=& \big\{ \psi \in Z^1(\Gmod; V); \; \exists v \in V \mbox{ such that
} \psi_\gamma = v|(1-\gamma) \big\}.
\end{eqnarray}
We give the arguments of the $\emph{cocycles}$
$\psi \in \Z^1(\Gmod;V)$ by a subscript. Furthermore we denote the
right action of $\Gmod$ on $V$ as $v \mapsto v | \gamma$ for
$\gamma \in \Gmod$ and $v \in V$. If the cohomology group is clear we
use the notation $[\ps]$ for the cohomology class of the
cocycle~$\ps$.

The \emph{first parabolic cohomology group} is the subgroup of
$H^1(\Gmod; V)$ given by
\begin{eqnarray}
\label{B4.2}
H_\parb^1(\Gmod; V)
&=& Z_\parb^1(\Gmod; V) / B^1(\Gmod; V),\\
\nonumber
Z_\parb^1(\Gmod; V)
&=& \big\{ \psi \in Z^1(\Gmod; V); \; \exists v \in V \mbox{ such that
} \psi_T = v|(1-T) \big\}.
\end{eqnarray}
If $W \supset V$ is a larger $\Gmod$-module, then the
\emph{first mixed parabolic cohomology
group} is given by
\begin{eqnarray}
\label{B4.3}
H_\parb^1(\Gmod; V,W)
&\!\!\!=\!\!\!& Z_\parb^1(\Gmod; V,W) / B^1(\Gmod; V),\\
\nonumber
Z_\parb^1(\Gmod; V,W)
&\!\!\!=\!\!\!& \big\{ \psi \in Z^1(\Gmod; V); \; \exists v \in W
\mbox{ such that } \psi_T = v|(1-T) \big\}.
\end{eqnarray}
Note that $H^1_\parb(\Gmod;V,W)$ is a subspace of $H^1(\Gmod;V)$, and
that there is a natural map $H^1_\parb(\Gmod;V,W) \longrightarrow
H^1_\parb(\Gmod;W)$.

\begin{remark}
\label{B4.4}
In the definitions above we use that $\Gmod \backslash \HH$ has only
one $\Gmod$-class of cusps represented by~$\infty$, and that the
subgroup $\Gmod_\infty\subset\Gmod$ fixing~$\infty$ is generated
by~$T$. We refer to \cite[Section~10]{BLZ} for a discussion of
parabolic cohomology for more general discrete subgroups of
$\PSL \R$.
\end{remark}

\begin{remark}Any cocycle is determined by its values on a set of
generators, so by $\ps_S$ and $\ps_T$ for~$\Gmod$. The relations
\eqref{B1.1} determine the relations
\begin{equation}\label{S-ST-rel}
\ps_S |(1+S) \;=\; 0\,,\qquad \ps_{T^{-1}S}|\left( 1 + T^{-1}S+
ST\right)\;=\;0\,.\end{equation}
For a parabolic cocycle $\ps$ we can arrange $\ps_T=0$ without
changing the cohomology class. The resulting cocycle is determined by
its value on $S$, subject to the relations
\begin{equation}\label{parcocrel}
\ps_S|_{2s}(1+S)\;=\;0\,,\qquad \ps_S \;=\; \ps_S|_{2s}(T+\TT)\,.
\end{equation}
\end{remark}

\subsection{Eigenfunctions of the Mayer operator and parabolic
cohomology}\label{sect-LM-coh}
In the Introduction we have already mentioned that eigenfunctions of
$\LLMayer s$ with eigenvalue $1$ or $-1$ give rise to elements of
$\FEd_s(\C')_\om$, where for $X\subset \C$:
\begin{align}\label{FEd-def}
\FEd_s(X)_\om &\;=\; \bigl\{ \text{analytic }P:X\rightarrow\C \;:\; P
= P|_{2s}T + P|_{2s}\TT\\
\nonumber
&\qquad \qquad \text{ on } X\cap T^{-1} X \cap (T')^{-1}X\bigr\}\,.
\end{align}
By analytic on $X\subset\R$ we mean \emph{real analytic}. For open
$X\subset \C$, analytic means \emph{holomorphic}.

The cocycle condition $\ps_S=\ps_S|(T+\nobreak\TT)$ in
\eqref{parcocrel} is similar to the three term equation \eqref{3te}.
In \cite[Section~13]{BLZ}, various aspects of the relation between
$\FEd_s(\C')_\om$ and cohomology are discussed. For the present paper
it is important that under assumption~\eqref{B2.2}:
\begin{equation}
H^1_\parb(\Gmod;\cV^\om_s,\cV^{\smp}_s) \;\cong\;
H^1_\parb(\Gmod;\cV^{\smp}_s) \;\cong\;
\FEd_s(0,\infty)_\om^{\mathrm{simple}}\,.
\end{equation}
The third superscript \emph{simple} indicates that we impose on
$P \in \FEd_s(0,\infty)_\om$ an asymptotic behavior at the end points
of $(0,\infty)$:
\begin{equation}\label{as-exp}
P(x) \;\sim\; \sum_{m=-1}^\infty c_m^\infty
x^{-2s-m}\quad(x\uparrow\infty)\,,\qquad P(x) \;\sim\;
\sum_{m=-1}^\infty c_m^0 x^m\quad(x\downarrow0)\,.
\end{equation}
This is a one-sided version of the behavior at singularities of
elements of~$\cV^{\smp}_s$ defined in \S\ref{B3}.

In \S\ref{sect-efMa} we shall show:
\begin{proposition}\label{prop-efLM-coh}Let $s\in \C$, $0<\re s<1$,
$s\neq\frac12$. The space $\ker(\LLMayer
s-\nobreak 1) \oplus \ker(\LLMayer s+ \nobreak 1)$ is in bijective
correspondence to $H^1_\parb(\Gmod;\cV^\om_s,\cV^{\smp}_s)$.
\end{proposition}

\subsection{Eigenfunctions of the nearest integer transfer operator
and hyperbolic cohomology}\label{sect-hypcoh}
The definitions \eqref{B4.2} and \eqref{B4.3} are related to the
$\Gmod$-orbit $\PPP\Q$ of cusps in $\PPP\R$. The element $T\in \Gmod$
generates the subgroup $\Gmod_\infty$ of $\Gmod$ fixing the element
$\infty$ in this orbit. Let us now work with what we would like to
call the \emph{Fibonacci orbit} $\cH = \Gmod (-\gr)\subset\PPP\R$,
where $\gr=\frac{1+\sqrt 5}2$ is the golden ratio. The hyperbolic
 element $TST^2\in \Gmod$ generates the subgroup $\Gmod_{-\gr}$
of~$\Gmod$ fixing $-\gr$.
\begin{definition}For $\Gmod$-modules $W\supset V$:
\begin{eqnarray}
\label{hypcoh}
H_\cH^1(\Gmod; V,W)
&\!\!\!=\!\!\!& Z_\cH^1(\Gmod; V,W) / B^1(\Gmod; V),\\
\nonumber
Z_\cH^1(\Gmod; V,W)
&\!\!\!=\!\!\!& \big\{ \psi \in Z^1(\Gmod; V): \; \exists v \in W
\;:\; \psi_{TST^2} = v|(1-TST^2) \big\},
\end{eqnarray}
and $H^1_\cH(\Gmod;V):=H^1_\cH(\Gmod;V,V)$.
\end{definition}

In particular, $H^1_\cH(\Gmod;\cV^{\fs}_s)$ is a subspace of
$H^1(\Gmod;\cV^{\fs}_s)$. The inclusion $\cV^\om_s \hookrightarrow 
\cV^{\fs}_s$ induces a linear map
$H^1(\Gmod;\cV^\om_s) \longrightarrow H^1(\Gmod;\cV^{\fs}_s)$.

Let us also define for $X\subset \C$:
\begin{align}
\FEv_s(X)_\om &\;=\; \biggl\{ \text{analytic }g:X\rightarrow\C \;:\;\\
\nonumber
&\qquad\qquad \hbox{} g+ g|_{2s}ST^2 \;=\; g|_{2s}T^{-1} +
g|_{2s}T^{-1}ST^{-2}\\
\nonumber
&\qquad\qquad \hbox{} \text{ on }X\cap T^{-2}SX \cap TX \cap T^2STX
\biggr\}\,,
\end{align}
with the same convention concerning analyticity as
in~\S\ref{sect-LM-coh}. We shall prove in
\S\ref{D2}:
\begin{theorem}\label{thm-th}Let $s\in \C$, $0<\re s<1$,
$s\neq\frac12$. There is an injective map
$\th : \FEv_s(-\gr^2,\gr) \longrightarrow
H^1_\cH(\Gmod;\cV^{\fs}_s)$. The image
$\th \left( \FEv_s(-\gr^2,\gr) \right) \subset
H^1(\Gmod;\cV^{\fs}_s)$ is equal to the image of
$H^1(\Gmod;\cV^\om_s)$ in $H^1(\Gmod;\cV^{\fs}_s)$.
\end{theorem}

\begin{proposition}\label{prop-an-coh}The natural map
$H^1(\Gmod;\cV^\om_s)\longrightarrow H^1(\Gmod;\cV^{\fs}_s)$ is
 injective.
\end{proposition}
\begin{proof}Let
$\ps\in Z^1_\parb(\Gmod;\cV^\om_s)$ such that $\ps_\g = f|(1-\g)$ for
all $\g\in\Gmod$ for some $f\in \cV^{\fs}_s$. {}From
$f|(1-\nobreak T)=\ps_T\in \cV^\om_s$ it follows that the set of
singularities $\sing f$ can contain at most the point~$\infty$;
otherwise $\sing f$ would be infinite. Hence
$\sing{f|S} \subset \{0\}$. {}From $f-f|S =\ps_S\in \cV^\om_s$ we
conclude that $f$ has no singularities at all, i.e.,
$f\in \cV^\om_s$. Hence $[\ps]=0$ in $H^1(\Gmod;\cV^\om_s)$.
\end{proof}
We now have the following system of injective maps:
\[ \xymatrix{ \FEv_s(-\gr^2,\gr) \ar@{^{(}->}[r]^\th
& H^1_\cH(\Gmod;\cV^{\fs}_s) \ar@{^{(}->}[r]
& H^1(\Gmod;\cV^{\fs}_s)\\
& H^1_\parb(\Gmod;\cV^\om_s,\cV^{\smp}_s) \ar@{^{(}->}[r]
& H^1(\Gmod;\cV^\om_s) \ar@{^{(}->}[u] }\]
\begin{theorem}\label{thm-par}Let $s\in \C$, $0<\re s<1$,
$s\neq\frac12$. The kernel of $\LLt s-1$ determines a subspace of
$\FEv_s(-\gr^2,\gr)_\om$ that is mapped by $\th$ onto the image of
$H^1_\parb(\Gmod;\cV^\om_s,\cV^{\smp}_s)$ in
$H^1(\Gmod;\cV^{\fs}_s)$.
\end{theorem}

This establishes a bijective map between $\ker (\LLt s-\nobreak 1)$
and the cohomology group $H^1_\parb(\Gmod;\cV^\om_s,\cV^{\smp}_s)$.
We prove this theorem in~\S\ref{sect-eif-coh}.

The results in Proposition \ref{prop-efLM-coh}, Theorem~\ref{thm-th},
Proposition \ref{prop-an-coh} and Theorem~\ref{thm-par} imply
Theorem~\ref{thm-main}. See also Table~\ref{overview}.
\begin{table}
\[ \xymatrix@R+3ex{ *+[F]{\ker(\LLMayer s-1)\oplus \ker(\LLMayer s+1)}
  \ar@(dl,ul)@/_/[d]_{\txt{\small \cite[Chap.~IV, \S3]{LZ01},
  Prop.~\ref{prop-Ma-FE}, \ \\\small restriction and extension \ }} \\
*+[F]{\FEd_s(0,\infty)_\om^{\mathrm{simple}} }
\ar@(dl,ul)@/_/[d]_{\txt{\small Prop.~\ref{prop-FE3-coh}, separation
  of singularities \ }} \ar@(ur,dr)@/_/[u]_{\txt{\small \
  Prop.~\ref{prop-Ma-FE}, bootstrap in\\\small
  \cite[Chap.~III,\S4]{LZ01}}} \\
*+[F]{H^1_\parb(\Gmod;\cV^\om_s,\cV^{\smp}_s)}
  \ar@{<->}[d]|{\txt{\small Prop.~\ref{prop-an-coh},
  $H^1(\Gmod;\cV^\om_s) \rightarrow H^1(\Gmod;\cV^{\fs})$ is
  injective}} \ar@(ur,dr)@/_/[u]_{\txt{\small \
  Prop.~\ref{prop-FE3-coh}, one-sided average $\avp$}} \\
*+[F]{\txt{Image of $H^1_\parb(\Gmod;\cV^\om_s,\cV^{\smp}_s)$ in
$H^1(\Gmod;\cV^{\fs}_s)$ }} \ar@(dl,ul)@/_/[d]_{\txt{\small
  Prop.~\ref{prop-thsurj}, one-sided average $\av_{TST^2}^+$, \ \\
  \small \S\ref{sect-eif-coh}, one-sided averages}} \\
*+[F]{\txt{subspace of $\FEv_s(-\gr^2,\gr)_\om$}}
  \ar@(dl,ul)@/_/[d]_{\txt{\small Lemma~\ref{lem-coc-FE}, cocycle
  relations, \ \\\small \S\ref{sect-eif-coh} }}
  \ar@(ur,dr)@/_/[u]_{\txt{\small \ Prop.~\ref{prop-om}, separation of
  singularities, \\\small \S\ref{sect-eif-coh}, one-sided averages}} \\
*+[F]{\txt{solutions of $\vec g|_{2s}\LL = \vec g$ on intervals} }
  \ar@(dl,ul)@/_/[d]_{\txt{\small Prop.~\ref{prop-bootstrap}, bootstrap
  \ }} \ar@(ur,dr)@/_/[u]_{\txt{\small \ Prop.~\ref{prop-thh}, cocycle
  relations, \\\small \S\ref{sect-eif-coh}}} \\
*+[F]{ \ker(\LLt s-1) } \ar@(ur,dr)@/_/[u]_{\txt{\small \
  Prop.~\ref{C3.2}, restriction and extension}} } \]
\caption{Overview of the steps in the proof of
Theorem~\ref{thm-main}.} \label{overview}
\end{table}

\subsection{Automorphic forms and cohomology groups} \label{B7}
In this note we work with transfer operators and cohomology groups. In
\cite{LZ01} and \cite{BLZ} the main theme is
the relation between period functions, automorphic forms and
cohomology. We mention the relevant facts as background material.

We denote by $\FEd_s(\C')^0_\om$ the subspace of
$P\in \FEd_s(\C')_\om $, as defined in \eqref{FEd-def}, with
 $\C'=\C\setminus(-\infty,0]$, that satisfy $P(x) = O(1)$ as
$x\downarrow0$ and $P(x)=O(x^{-2s})$ as
$x\rightarrow\infty$. The main theorem in \cite{LZ01} states that the
space $\FEd_s(\C')_\om^0 \cong H^1_\parb(\Gmod;\cV^{\sa}_s)$ is in
bijective correspondence with the space of Maass cusp forms with
spectral parameter~$s$. A Maass cusp form is a function
$u:\HH\longrightarrow\C$ satisfying $u(\g z) = u(z)$ for all
$\g\in \Gmod$ that is given by a convergent Fourier expansion
\begin{equation}
u(x+iy) \;=\; \sum_{n\neq 0} A_n e^{2\pi i nx} \sqrt y\,
K_{s-1/2}(2\pi|n|y)\,.
\end{equation}
The space $M_s^0$ of such Maass cusp forms is known to be non-zero
only for a discrete set of values of $s$ satisfying $\re s=\frac12$,
$s\neq\frac12$.

A slightly larger space of $\Gmod$-invariant functions is $M_s^1$,
consisting of the $\Gmod$-in\-vari\-ant $u$ on~$\HH$ with a
converging Fourier expansion
\begin{equation}
u(x+iy) \;=\; A_0 y^{1-s} + \sum_{n\neq 0} A_n e^{2\pi i nx} \sqrt y\,
K_{s-1/2}(2\pi|n|y)\,.
\end{equation}
This space is equal to $M_s^0$ for $\re s=\frac12$, $s\neq \frac12$.
For values of~$s$ with $0<\re s<\frac12$ such that $\zeta(2s) = 0$
the residue of the Eisenstein series is an element of $M_s^1$.
The results in \cite[Section~11]{BLZ} show that for $0<\re s<1$,
$s\neq\frac12$ there is a bijective correspondence between $M_s^1$
and $H^1_\parb(\Gmod;\cV^{\smp}_s)$.
These spaces are finite dimensional, and zero for most values of~$s$.
All elements of $M_s^1$ are eigenfunctions of the hyperbolic Laplace
operator: $-y^2\left( \partial_y^2+\partial_x^2\right)u = s
(1-\nobreak s)u$.

The conclusion is that the eigenfunctions of $\LLMayer s$ with
eigenvalues $1$ and $-1$, and the eigenfunctions of $\LLt s$ with
eigenvalue~$1$ are in bijective correspondence to elements of the
space $M_s^1$. This gives a confirmation of the relation between
eigenfunctions of transfer operators and automorphic forms that we
know already from the relation via the Selberg zeta function. (See
\cite{Ma} and \eqref{B8.1}.)

%%%%%%%%%%%% Section 3 %%%%%%%%%%%%%%%

\section{The transfer operators on disks and intervals} \label{C}

In the context of dynamical systems one usually considers transfer
operators in Banach spaces of holomorphic functions on a disk. For
the relation to cohomology groups with values in principal series
spaces, it is more natural to consider the corresponding operators on
functions on intervals in $\R$ or $\PPP\R$. We discuss this relation
in \S\ref{sect-efMa} for the Mayer operator and in
\S\ref{sect-ef-nito} for $\LLt s$. In \S\ref{sect-4te} we derive the
four term equation~\eqref{4te} from the transfer operator $\LLt s$.

We start with a discussion of one-sided averages.

\subsection{One-sided averages} \label{C1}

Both in the definition of $\LLMayer s$ in \eqref{A.1} and in that of
$\LLt s$ in \eqref{A.5}, one recognizes infinite sums of the type
$f| \g T^n$ over infinitely many~$n\in \Z$ for a fixed $\g\in\PGL\Z$.
The one-sided averages
\begin{equation}
\label{C1.1}
\avp = \sum_{n=0}^\infty T^n \qquad \mbox{and} \qquad \avm =
-\sum_{n=-\infty}^{-1} T^n
\end{equation}
play also an important role in~\cite{BLZ}. In this subsection we
recall the relevant results.\medskip

Consider a function of the form $f=h|S$, where $h$ is holomorphic on a
neighborhood of~$0$. Then
\begin{eqnarray}\label{avdef}
\av_T^+(f) \;=\; f| \avp (z) &=& \sum_{n=0}^\infty ((z+n)^2)^{-s} \, h
\left( \frac{-1}{z+n} \right)
\quad\mbox{and} \\
\nonumber
\av_T^-(f) \;=\; f| \avm (z) &=& -\sum_{n=1}^\infty ((z-n)^2)^{-s} \,
h \left( \frac{-1}{z-n} \right)
\end{eqnarray}
converge absolutely if $\re s>\frac12$, and define $f|\avp$ as a
holomorphic function on a right half-plane, and $f|\avm$ on a left
half plane. If $h(0)=0$, the convergence is absolute for $\re s>0$.
Using the Hurwitz zeta function for the contribution of the constant
term of $h$ at~$0$, we obtain in general a meromorphic continuation,
with at most a first order singularity at $s=\frac12$ on $\re s>0$.
In this note we will understand $f|\avp$ and $f|\avm $ always in this
regularized sense. We have given two notations in~\eqref{avdef}. With
$f|\avp$ we stress that $\avp$ is an element of the completion of the
group ring of~$\Gmod$, for which we have made sense of the action on
certain functions by regularization. With $\avp(f)$ we emphasize that
this one-sided average defines an operator on suitable spaces of
functions. In this note we will use $f|\avp$ and $f|\avm$.

These one-sided averages satisfy
\begin{alignat}2 \label{avTrel1}
f|\avp |(1-T) &\;=\; f\,,&\qquad f|\avm |(1-T) &\;=\; f\,,\\
\label{avTrel2}
f|(1-T)|\avp &\;=\; f\,,&\qquad f|(1-T)|\avm
&\;=\; f\,,\\
\label{avTcom}
f|T|\avp &\;=\; f|\avp|T \,,&\qquad f|T|\avm&\;=\; f|\avm|T\,,\\
\label{avTsplit}
f + f|T |\avp &\;=\; f|\avp\,,&\qquad
-f|T^{-1} + f|T^{-1}|\avm&\;=\; f|\avm\,,
\end{alignat}
on suitable right half-planes, respectively left half-planes. These
relations hold trivially in the domain $\re s>\frac12$ of absolute
convergence, and survive under meromorphic continuation.

In particular, we consider these one-sided averages for
$f\in\cV^\om(I)$ where $I\subset\PPP\R$ is a neighborhood
of~$\infty$. Then $f$ has the form indicated above. Let us consider a
cyclic interval $I=(a,b)_c$ in $\PPP\R$ containing $\infty$. (This
means that $a>b$ in $\R$ and $(a,b)_c = (a,\infty)
\cup \{\infty\}\cup(-\infty,b)$.\@)
 As in \cite[Section~3]{BLZ} we have:
\begin{lemma}
\label{C1.2}Let $ I \supset (a,b)_c\ni \infty$. If $f\in
\cV^\om_s\bigl(I\bigr)$, then $f|\avp \in \cV^\om_s(a,\infty)$ and is
represented by a function holomorphic on a neighborhood of
$(a,\infty)$ containing a right half-plane, and
$f|\avm\in \cV^\om_s (\infty,b+\nobreak1)$ is represented by a
function holomorphic on a neighborhood of $(\infty,b+\nobreak1)$
containing a left half-plane.

There are constants $C_m^\ast$ for $m=-1,0,1,\ldots$ such that
\begin{alignat}2 \label{avp-as}
f|_{2s}\avp(x) &\;\sim\; \sum_{m=-1}^\infty C_m^\ast
x^{-m-2s}&&(x\uparrow\infty)\,,\\
\label{avm-as}
f|_{2s}\avm(x)&\;\sim\; \sum_{m=-1}^\infty C_m^\ast
x^{-m}|x|^{-2s}&\quad&
(x\downarrow-\infty)\,.
\end{alignat}
\end{lemma}
In particular, if $f\in \cV^\om_s=\cV^\om_s(\PPP\R)$, then
$f|\avpm \in \cV^\om_s(\R)$.\smallskip

The asymptotic behavior in \eqref{avp-as} and \eqref{avm-as} is
related to the singularity behavior \eqref{simple-sing} in the
definition of~$\cV^\smp_s$:
\begin{lemma}\textrm{\cite[Section~8]{BLZ}
}\label{8.17} For $f\in \cV^\om_s $ the following statements are
equivalent:
\begin{enumerate}
\item[i)] $f\in \cV^{\smp}_s|(1-T)$.
\item[ii)] $f|\avp=f|\avm$.
\end{enumerate}
\end{lemma}
\begin{proof}If $f|\avp=f|\avm$, then
$f|\av^\pm|S(z) \sim \sum_{m=-1}^\infty C_m^\ast x^m$ as
$\mp x \downarrow0$. Hence $f\in \cV^{\smp}_s$.

Conversely, suppose that $f=h|(1-\nobreak T)$ with
$h\in \cV^{\smp}_s$. Then $p_+=h-f|\avp$ and $p_-=h-f|\avm$ satisfy
$p_\pm |T=p_\pm$. Moreover, $p_+(x)$ has an asymptotic expansion of
the form \eqref{avp-as} as $x\uparrow \infty$. Hence the periodic
function $p_+$ vanishes. For $p_-$ let $x\downarrow-\infty$.
\end{proof}

\begin{lemma}
\label{C1.8}
Suppose that $b,c \in \cV^\omega_s$ satisfy
$ b | \avp + c | \avm \in \cV^\om_s$. Then
$b|\avp=b|\avm \in \cV^{\smp}_s$, and
$ c|\avp=c|\avm \in \cV^{\smp}_s$.
\end{lemma}
\begin{proof}
Relations \eqref{avTrel1} and \eqref{avTrel2} imply that
$p=c|\avp - c|\avm$ satisfies $p|_{2s}T=p$, hence $p$ is a periodic
function on~$\R$. Put $a= b | \avp + c | \avm$. As $x\uparrow\infty$,
the term $c|\avp$ has an asymptotic expansion as in \eqref{avp-as},
and $c|\avm= a-b|\avp$ also has an expansion of this type. (For $a$
we know that $a(x) = |x|^{-2s} $(analytic in $-1/x$).) Hence
$p(x) \sim x^{-2s} \left( p_{-1} x + p_0 + \cdots\right)$ as
$x\rightarrow\infty$. The periodicity implies that $p$ is bounded,
hence $p_{-1}=0$. Next $p(x) = \OO{x^{-2\re s}}$ implies $p=0$, hence
$c|\avp=c|\avm$. Now $c|\avpm \in \cV^{\smp}_s$, and also
$b|\avpm\in \cV^{\smp}_s$.
\end{proof}
\medskip

We can build one-sided averages for other elements of~$\Gmod$. If
$\eta\in \Gmod$ is hyperbolic, for instance $\eta=TST^2$, then the
averages $\av_\eta^+=\sum_{n\geq 0}\eta^n$ and
$\av_\eta^-=-\sum_{n\leq -1}\eta^n$ have the properties corresponding
to \eqref{avTrel1}--\eqref{avTsplit}, if they converge,
\cite[Section~7]{BLZ}. If the attracting fixed point $\om(\eta)$
of~$\eta$ is in the cyclic interval $I\subset\PPP\R$, then
$f|_{2s}\av_\eta^+$ converges without regularization for $\re s>0$
for all $f\in \cV^\om_s(I)$, and provides us with
$f|_{2s}\av_\eta^+ \in \cV^\om_s\left( I
\setminus\{\al(\eta)\}\right)$, where $\al(\eta)$ is the repelling
fixed point of~$\eta$. In particular,
$\av_\eta^+ : \cV^\om_s \rightarrow \cV^\om_s\left(\PPP\R\setminus\{\al(\eta)
\}\right)$, and similarly $\av_\eta^-:\cV^\om_s \rightarrow
\cV^\om_s\left(\PPP\R\setminus \{\om(\eta)\}\right)$.

\subsection{Eigenfunctions of the Mayer operator}\label{sect-efMa}
Eigenfunctions of $\LLMayer s$ with eigenvalue $\pm1$, or briefly
$(\pm1)$-eigenfunctions of $\LLMayer s$, can be related to
eigenfunctions of a similar operator on $\cV^\om_s(0,\infty)$. This
statement is almost contained in the results in \cite{LZ01}.
Nevertheless, we recall the main steps in the proof, since we want to
describe the correspondence between $\pm1$-eigenfunctions of
$\LLMayer s$ and $1$-eigenfunctions of $\LLt s$ explicitly.

In this subsection we also indicate how to prove
Proposition~\ref{prop-efLM-coh}, on the basis of results
in~\cite{BLZ}.
\medskip

Under the step $P=f|T$, the function
 $f\in \ker\left( \LLMayer s\mp1\right)$ corresponds to $P$
holomorphic on $|z-\nobreak 2|< \frac32$, and continuous on
$|z-\nobreak 2|\leq \frac 32$ satisfying $P|C\TT\avp = \pm P$, with
use of the notation introduced in \eqref{B1.3}. This implies on
suitable non-empty domains:
\begin{alignat*}2 P&\;=\; P|T \pm P|TC \,,&\qquad P|C&\;=\; \pm P\,,\\
P&\;=\; P|T + P|\TT\,.
\end{alignat*}
The last equality is the three term equation~\eqref{3te}. The
proposition in \cite[Chap.~IV, \S3]{LZ01} shows that these functions
extend to $\C'=\C\setminus(-\infty,0]$ and satisfy
$P(x) \sim c_{-1} x^{1-2s} + O(x^{-2\re s})$ as $x\uparrow\infty$ for
some $c_{-1}$. Since $P$ satisfies $P|T'\avp = P$ on $\C'$ and
$P|C=\pm P$, we have the asymptotic behavior \eqref{as-exp} near both
end points of $(0,\infty)$, with $c^0_m = \pm c^\infty_m$. Thus we
have $P\in \FEd_s^\pm(\C')_\om^{\mathrm{simple}}$, where the upper
index $\pm$ indicates the $(\pm1)$-eigenspace of $C$ in
$\FEd_s(\C')_\om$, and where the superscript \emph{simple} indicates
the subspace satisfying \eqref{as-exp}.

Conversely, starting with
$P\in \FEd_s^\pm(\C')_\om^{\mathrm{simple}}$, we have
$P|\TT\in \cV^\om_s\left( (0,-1)_c\right)$, where
$(0,-1)_c=(0,\infty)
\cup\{\infty\}\cup(-\infty,0)$ denotes a cyclic interval in $\PPP\R$.
Hence $P|\TT\avp\in \cV^\om_s(0,\infty)$. By~\eqref{avTrel1}:
\[ \left( P - P|\TT\avp\right)|(1-T) \;=\; P-P|T - P|\TT \;=\;0\,.\]
The asymptotic behavior of $P - P|\TT\avp$ near~$\infty$ shows that
this periodic function vanishes. So $P$ satisfies $P|\TT\avp=P$, or
 with use of the parity condition, $P|C\TT \avp=\pm P$. The function
$P|\TT\avp$ is holomorphic on a right half plane. With the parity
condition $P|C=\pm P$ this implies that $P$ is holomorphic on a wedge
of the form $|\arg z|<\e$. This suffices as the point of departure
for the second stage of the bootstrap procedure in \cite[Chap.~III,
\S4]{LZ01}, which gives a holomorphic extension of $P$ to $\C'$,
still satisfying $P|C\TT\avp = \pm P$. This leads to
$f=P|T \in \ker(\LLMayer s\mp\nobreak1)$.

Thus, we have an explicit bijective correspondence between the
following spaces:
\begin{align*}
&\ker \left( \LLMayer s\mp 1\right)\,,\\
& \ker\left(T'\avp -1:\cV^\om_s(0,\infty) \longrightarrow
\cV^\om_s(0,\infty)
\right) \\
\nonumber
&\qquad\qquad\hbox{} \cap \ker \left( C\mp1:\cV^\om_s(0,\infty)
\longrightarrow \cV^\om_s(0,\infty)
\right)\,,
\\
&\FEd_s^\pm(\C')_\om^{\mathrm{simple}}\,, \text{ and }
\FEd_s^\pm(0,\infty)_\om^{\mathrm{simple}}\,.
\end{align*}
As in \cite[Chap.~I, \S3]{LZ01}, we have
$\FEd_s(\C')_\om^{\mathrm{simple}} =
\FEd_s^+(\C')_\om^{\mathrm{simple}} \oplus \FEd_s^-(\C')_\om^{\mathrm{simple}}$.
This gives the following result:
\begin{proposition}\label{prop-Ma-FE}Let $0<\re s<1$, $s\neq\frac12$.
There is an explicit bijection between the following two spaces:
\begin{align*}
&\ker(\LLMayer s-\nobreak1) \oplus \ker(\LLMayer s+\nobreak1)\,,\\
&\FEd_s(0,\infty)_\om^{\mathrm{simple}}\,.
\end{align*}
\end{proposition}

To complete the proof of Proposition~\ref{prop-efLM-coh} we have to
establish a relation between $\FEd_s(0,\infty)^{\mathrm{simple}}_\om$
and $H^1_\parb(\Gmod;\cV^\om_s,\cV^{\smp}_s)$. The least explicit
step in the proof is provided by the following result:
\begin{proposition}\label{D2.6} If $f \in \cV^\fs_s$ satisfies
$\sing{f} \subset \{\xi,\eta\}$ for two different points $\xi$ and
$ \eta$ in $ \PPP\R$, then there are $f_\xi, f_\eta \in \cV^\fs_s$
such that $f = f_\eta- f_\xi$ and $\sing{f_\xi} \subset \{\xi\}$,
$\sing{f_\eta} \subset \{\eta\}$. The functions $f_\xi$ and $f_\eta$
are not unique. The freedom consists of adding the same element of
$\cV^\omega_s$ to both functions.
\end{proposition}
\begin{proof}[Sketch of a proof] This follows from, e.g.,
\cite[Theorem~1.4.5]{Ho}. See \cite[Section~12]{BLZ} for the
application to elements of~$\cV^{\fs}_s$.

The idea is to use another model of the principal series, in which the
elements of $\cV^\om_s$ correspond to functions holomorphic on an
annulus in~$\C$ containing the unit circle. The function $f$ in the
proposition is represented by a holomorphic function on an open set
$\Om\subset\C$ containing the unit circle minus the points $\tilde\x$
and $\tilde\eta$ corresponding to $\x$ and~$\eta$. Write
$\Om = \Om_1\cap \Om_2$ with $\tilde \eta\in \Om_1$,
$\tilde \x\in \Om_2$. Apply \cite[Theorem~1.4.5]{Ho} with
$g_{1,2} = f$ to obtain $f=g_1-g_2$ on~$\Om$ with $g_j$ holomorphic
on~$\Om_j$.
\end{proof}

\begin{proposition}\label{prop-FE3-coh}Let $0<\re s<1$,
$s\neq \frac12$. There is an explicit bijection between
$\FEd_s(0,\infty)^{\mathrm{simple}}_\om$ and
$H^1_\parb(\Gmod;\cV^\om_s,\cV^{\smp}_s)$.
\end{proposition}
\begin{proof}Suppose that
$P\in \FEd_s(0,\infty)^{\mathrm{simple}}_\om$. We extend it to
$\tilde P\in \cV^{\fs}_s$ by
\begin{equation} \tilde P\;=\;P\text{ on }(0,\infty)\,,\qquad \tilde
P\;=\;-P|S\text{ on }(-\infty,0)\,.
\end{equation}
 So $\sing {\tilde P} \subset \{0,\infty\}$. By separate computations
 on $(-\infty,-1)$ and $(-1,0)$ we conclude that $\tilde P$ satisfies
 \eqref{3te} on $\PPP\R\setminus\{\infty,-1,0\}$. Thus, we obtain a
 cocycle $\ps \in Z^1_\parb(\Gmod;\cV^{\fs}_s)$, determined by
\[ \ps_T \;=\;0\,,\qquad \ps_S\;=\; \tilde P\,.\]

Since
$\tilde P|\TT \in \cV^\om_s\left(\PPP\R\setminus\{-1,0\}\right)$, we
have $\tilde P |\TT \avp = P |\TT\avp \in \cV^\om_s(0,\infty)$ and
$\tilde P|\TT\avm \in \cV^\om_s(-\infty,0)$. We have indicated
earlier in this subsection that the asymptotic behavior of $P(x)$ as
$x\uparrow\infty$ implies that $\tilde P|\TT\avp = \tilde P$ on
$(0,\infty)$. The asymptotic behavior of $P(x)$ as $x\downarrow 0$
implies a similar asymptotic behavior of $\tilde P(x)$ as
$x\downarrow-\infty$, which in turn implies analogously that
$\tilde P |\TT\avm = \tilde P$ on $(-\infty,0)$. Consulting
\eqref{avp-as} and \eqref{avm-as} we conclude that $\tilde P$ has the
same coefficients in its expansions for both directions of approach
to $\infty\in \PPP\R$. The fact that $\tilde P|S=-\tilde P$ implies
the same statement at~$0$. Hence $\tilde P \in \cV^{\smp}_s$ and
$\ps\in Z^1_\parb(\Gmod;\cV^{\smp}_s)$. It is the unique cocycle in
its cohomology class in $H^1_\parb(\Gmod;\cV^{\smp}_s)$, since
$(\cV^{\smp}_s)^T=\{0\}$, as is shown in \cite[Section~8]{BLZ}.

Proposition~\ref{D2.6} implies that there are
$F_\infty, F_0\in \cV^{\fs}_s$ such that $\tilde P=F_\infty-F_0$ with
$\sing {F_\infty}\subset\{\infty\}$, $\sing{F_0}\subset\{0\}$. Since
$F_\infty$ has the same type of asymptotic behavior at~$\infty$ as
$\tilde P$, we conclude that $F_\infty, F_0\in \cV^{\smp}_s$.
Moreover, $\tilde P|S=-\tilde P$ implies that $F_0 = F_\infty+\al$
for some $\al \in \cV^\om_s$. {}From $\sing{F_\infty|\TT}\subset
\{-1\}$ and $\sing{F_0|S\TT}\subset\{0\}$, it follows that
\begin{align*}
\sing{\tilde P|(1-T) }&\;=\; \sing{ F_\infty|\TT
-F_\infty|S\TT-\al|\TT}\; \subset\; \{-1,0\}\,,\\
\sing{ F_\infty|(1-T)} &\;=\;\sing{ \tilde P|(1-T) + F_0|(1-T)}
\;\subset\; \{-1,0\}\,.
\end{align*}
On the other hand $\sing{F_\infty}\subset\{\infty\}$ implies that
$\sing{F_\infty|(1-T)}\subset\{\infty\}$. The conclusion is that
$F_\infty|(1-T)
\in \cV^\om_s$. We conclude that the cocycle
\[\tilde\ps: \g\mapsto \ps_\g-F_\infty|(1-\g)\]
takes values in $\cV^\om_s$. The freedom in the choice of $F_\infty$
and $F_0$ amounts to the freedom of choosing $\tilde\ps$ in its class
in $H^1_\parb(\Gmod;\cV^\om_s,\cV^{\smp}_s)$.

Conversely, we start with a cocycle $\tilde \ps \in
Z^1_\parb(\Gmod;\cV^\om_s,\cV^{\smp}_s)$. By the definition of mixed
parabolic cohomology in \eqref{B4.3}, there exists
$v\in \cV^{\smp}_s$ such that $\tilde \ps_T = v|(1-\nobreak T)$. A
possible choice is $\tilde\ps_T|\avp$, which coincides with
$\tilde\ps_T|\avm$ according to Lemma~\ref{8.17}, since
$\ps_T\in \cV^{\smp}_s$. Lemma~\ref{C1.8} shows that
$v=\tilde\ps_T|\avp$ is the sole possibility.

Now $\ps_\g = \tilde\ps_\g - \tilde\ps_T|\avp|(1-\nobreak\g)$
determines the unique cocycle in the cohomology class of $\tilde\ps$
in $H^1_\parb(\Gmod;\cV^{\smp}_s)$
vanishing on~$T$. Note that $\ps$ does not change if we change
$\tilde\ps$ in its cohomology class in
$H^1_\parb(\Gmod;\cV^\om_s,\cV^{\smp}_s)$.
We have $\ps_T=0$ and
$\ps_S = \tilde\ps_S -\tilde\ps_T|\avp|(1-\nobreak S) \in
\cV^{\smp}_s$ with singularities contained in $\{0,\infty\}$.
Restriction of $\ps_S$ to $(0,\infty)$ gives an element of
$\FEd_s(0,\infty)^{\mathrm{simple}}_s$.

Noting that $F_\infty = - \tilde\ps_T |\avp$, we check that both
procedures are inverse to each other.
\end{proof}
There might be cohomology classes in $H^1_\parb(\Gmod;\cV^{\smp}_s)$
that are not represented by a cocycle $\ps$ such that
$\sing{\ps_S}\subset \{0,\infty\}$. In \cite[Section 12]{BLZ} it
takes work to show that such classes do not exist.

\subsection{Eigenfunctions of the nearest integer transfer
operator}\label{sect-ef-nito}
In this subsection we relate $1$-eigenfunctions of $\LLt s$ in
\eqref{A.5} to vectors of real analytic functions on intervals. We
note that $\LLt s$ is given by
\begin{equation}\label{LL-f}
\vec g \mapsto \vec g| \LL\,,\quad\text{ where } \quad \LL\;=\;
\rMatrix{ST^3\avp}{ST^2\avp}{-ST^{-1}\avm}{-ST^{-2}\avm}\,,
\end{equation}
for $\vec g=(g_1,g_2)$ continuous on
$D=\left\{z\in \C\;:\; |z|\leq 1\right\}$ and holomorphic on the
interior~$\mathring D$. The transition from the operator notation in
\eqref{A.5} to the right module notation here causes a transition
from column vectors to row vectors and a transposition of the matrix.

For the first component $g_1$ of $\vec g =(g_1,g_2)$ in the domain of
$\LL$ the translate $g_1|ST^3$ is holomorphic on the region
$|z+\nobreak 3|> 1$ in $\PPP\C$, hence $g_1|ST^3\avp$ is at least
defined as a holomorphic function on the right half plane $\re z>-2$
in~$\C$. Proceeding similarly with the other components we find that
$\vec g|\LL$ is holomorphic on $-1<\re z<1$, which contains the
interior of the unit disk~$D$. Considering two terms in the infinite
sum separately, with \eqref{avTsplit}, we see that $\vec g|\LL$ is
continuous on the unit disk $D$, its boundary included. Restriction
of $\vec g$ to a neighborhood of $(-1,1)$ provides us with a solution
of $\vec g|\LL = \vec g$ in the vector valued analytic functions
on~$(-1,1)$.

We shall see that the components of this restriction to $(-1,1)$ can
be extended to a larger interval, and that conversely solutions of
$\vec g = \vec g|\LL$ on the resulting interval come from
$1$-eigenfunctions of~$\LLt s$.

It turns out that a crucial role is played by the intervals
$(-\gr^2,\gr)$ and $(-\gr,\gr^2)$, where $\gr$ denotes the golden
ratio $\gr=\frac{1+\sqrt 5}2$.
\begin{proposition}
\label{C3.2}
Let $a_1,b_2 \in \left( \gr^{-2},\gr^2 \right)$ and
$a_2,b_1\in (0,1]$. Suppose that $f_1\in \cV^\om_s(-a_1,b_1)$ and
$f_2\in \cV^\om_s(-a_2,b_2)$ satisfy
$\left( f_1,\, f_2\right)|_{2s}\LL = \left( f_1,\, f_2 \right)$. Then
$f_1 $ is the restriction of $h_1 \in \cV^\om_s(-\gr^2,\gr)$ and
$f_2$ of $h_2\in \cV^\om_s(-\gr,\gr^2)$ such that $(
h_1, h_2)|\LL = ( h_1, h_2)$, and $h_1=h_2|T$. The values of
$h_1(x)$, respectively $h_2(x)$ for each given $x\in
(-\gr^2,\gr)$, respectively $x\in (-\gr,\gr^2)$, can be expressed in
values of $f_1$ and $f_2$ by a finite number of applications of the
relation $(f_1,f_2)=(f_1,f_2)|_{2s}\LL$.
\end{proposition}
\begin{proof}
By analyticity the eigenfunction equation extends from given open
intervals to larger ones, and relation \eqref{avTcom} implies
$h_1=h_2|T$. The statement that requires work is the extension of the
domains.

Denote $(f_1,f_2)|\LL$ by $(\tilde f_1, \tilde f_2) $. We have
\begin{alignat*}2 f_1|ST^3 &\;\in\; \cV^\om_s\bigl(
(\txtfrac1{a_1}-3,\txtfrac{-1}{b_1}-3)_c\bigr)\,,&\quad
 f_2|ST^{-1}&\;\in\;\cV^\om_s \bigl(
 (\txtfrac1{a_2}+1,\txtfrac{-1}{b_2}+1)_c\bigr)\,,\\
f_1|ST^2 &\;\in\; \cV^\om_s\bigl(
(\txtfrac1{a_1}-2,\txtfrac{-1}{b_1}-2)_c\bigr)\,,&\quad
 f_2|ST^{-2}&\;\in\;\cV^\om_s \bigl(
 (\txtfrac1{a_2}+2,\txtfrac{-1}{b_2}+2)_c\bigr)\,,
\end{alignat*}
where $(x,y)_c=(x,\infty)\cup\{\infty\}\cup(-\infty,y)$ is the
notation for cyclic intervals in $\PPP\R$. With Lemma~\ref{C1.2}:
\begin{alignat*}2 f_1|ST^3\avp &\;\in\; \cV^\om_s\bigl(
\txtfrac1{a_1}-3,\infty\bigr)\,,&\quad f_2|ST^{-1}\avm&\;\in\;
\cV^\om_s\bigl( -\infty, \txtfrac{-1}{b_2}+2\bigr)\,,\\
f_1|ST^2\avp &\;\in\; \cV^\om_s\bigl(
\txtfrac1{a_1}-2,\infty\bigr)\,,&\quad f_2|ST^{-2}\avm&\;\in\;
\cV^\om_s\bigl( -\infty, \txtfrac{-1}{b_2}+3\bigr)\,,\end{alignat*}
which implies
\begin{alignat*}2 \tilde f_1 \;\in\; \cV^\om_s&\bigl(
\txtfrac1{a_1}-3, \txtfrac{-1}{b_2}+2\bigr)\,,&\quad \tilde f_2
\;\in\; \cV^\om_s&\bigl( \txtfrac1{a_1}-2,
\txtfrac{-1}{b_2}+3\bigr)\,.
\end{alignat*}
The end points of the domains are transformed according to $-\tilde
a_1=\frac1{a_1}-3$, $\tilde b_1 = \tilde b_2-1$,
$-\tilde a_2= -\tilde a_1+1$, $\tilde b_2 = 3-\frac1{b_2}$. Iterating
this, the $a_1$ starting in $\left( \gr^{-2},\gr^2\right)$ form a
sequence increasing to $\gr^2$, and similarly for the $b_2$. This
leads to the extension indicated in the proposition.
\end{proof}
This proposition shows that $1$-eigenfunctions of $\LLt s$ extend to
vectors of the form $(g,g|T)$ with $g$ holomorphic on a neighborhood
of $(-\gr^2,\gr)$, such that the relation $(g,g|T) |\LL = (g,g|T)$ is
valid on a neighborhood of $(-\gr,\gr)$.

\begin{proposition}\label{prop-bootstrap}Suppose that
$g \in \cV^\om_s$ satisfies
\begin{equation}
\label{1ef}
(g, g|T^{-1})|\LL = (g,g|T^{-1})
\end{equation}
on $(-\gr,\gr)$. Then $g$ extends holomorphically to a neighborhood of
the closed unit disk $D$ and \eqref{1ef} holds on that neighborhood.
\end{proposition}
\begin{proof}Since $g$ is real analytic on the interval $(-\gr^2,\gr)
\supset [-1,1]$, there is a complex $\e$-neighborhood $U$ of $[-1,1]$
to which $g$ extends as a holomorphic function. Relation \eqref{1ef}
stays valid on this neighborhood.

Denote $\vec g = (g,g|T^{-1})$. We have
\begin{equation}\label{A-eq}
\vec g|\LL \;=\; \sum_{n\geq 0} \vec g | A_n \,,\qquad A_n \;=\;
\rMatrix {ST^{3+n}} {ST^{2+n}} {ST^{-2-n}} {ST^{-3-n}}\,.
\end{equation}
For sufficiently large $n$ the four images $(A_n)_{i,j} D$,
$i,j\in \{1,2\}$, are contained in the given neighborhood $U$ of
$[-1,1]$. Lemma~\ref{C1.2} and repeated application of relation
\eqref{avTsplit} show that there is a tail of the series in
 \eqref{A-eq} representing a holomorphic function on a neighborhood of
$D$. Thus we have on an open neighborhood $\Om$ of $(-\gr,\gr)$:
\[ \vec g|\LL \;=\; \sum_{n=0}^N g|A_n + \left( \text{holomorphic on a
neighborhood of $D$}\right)\,.\]

For $g$ in the remaining terms $g|A_n$ we substitute \eqref{A-eq}
again. Repeating this process, we obtain for each $k\geq 1$ on an
open neighborhood $\Om_k$ of $(-\gr,\gr)$:
\begin{align}\label{L-iter}
 \vec g(z) &\;=\; \sum_{n_1=0}^{N_1} \cdots \sum_{n_k=0}^{N_k} g|
A_{n_k}\cdots A_{n_1}(z)\\
\nonumber
&\qquad\hbox{}
 + \left( \text{holomorphic on a neighborhood of $D$}\right)\,.
\end{align}
The neighborhood $\Om_k$ is increasing in~$k$.

The matrix elements of $A_{n_k}\cdots A_{n_1}$ are of the form
$ST^{a_k}ST^{a_{k-1} }\cdots S T^{a_1}$, where $a_j\in \Z$,
$|a_j|\geq 2$, and $a_ja_{j+1}=-4$ if $|a_j|=|a_{j+1}|=2$. (The
transfer operator $\LLt s$ is designed to reflect this condition on
the $a_j$ in nearest integer continuous fraction expansions.)
Each $ST^{a_j}$ maps the unit disk $ D$ into itself. For the imaginary
part $y$ of $z\in D$ we have
\[ \im {S T^{a_j} z} \;=\; \frac y{|z+a_j|^2} \;\leq \; \frac
y{(|a_j|-1)^2}\,.\]
So each $ST^{a_j}$ with $|a_j|\geq 3$ decreases the imaginary part
with at least a factor~$4$. If $a_j=\pm 2$ the imaginary part does
not increase and if moreover $j<k$, we know that $a_{j+1}\neq a_j$.
One checks for $a_j a_{j+1}\geq 6$ and $a_j a_{j+1} \leq -4$
separately that the imaginary part decreases at least with a factor
$4$ under $ST^{a_{j+1}}ST^{a_j}$. So if $k$ is sufficiently large,
all $ST^{a_k}\cdots ST^{a_1}z$ with $z\in D$ are contained in the
$\e$-neighborhood $U$ of $[-1,1]$ we started with. For such $k$, all
explicit terms in \eqref{L-iter} can be absorbed in the last
term.\end{proof}

These propositions imply:
\begin{corollary}Let $0<\re s<1$, $s\neq\frac12$. Restriction and
analytic extension give a bijective correspondence between
$\ker(\LLt s - \nobreak1)$ and the space of solutions in
$\cV^\om_s(-\gr^2,\gr)
\times \cV^\om_s(-\gr,\gr^2)$ of
\begin{equation}
(g_1,g_2)|_{2s} \LL \;=\; (g_1,g_2)\,.
\end{equation}
\end{corollary}

\subsection{Four term equation}\label{sect-4te}

\begin{proposition}\label{prop-1ei-4te}The solutions in
$\cV^\om_s(-\gr^2,\gr) \times \cV^\om_s(-\gr,\gr^2)$ of
$\vec g = \vec g |_{2s}\LL$ are of the form $(g,g|T^{-1})$ with
$g\in \FEv_s(-\gr^2,\gr)_\om$.
\end{proposition}
\begin{proof}The eigenfunction relations
\begin{alignat*}2 h_1 &\;=\; & h_1|ST^3\avp &-h_2|ST^{-1}\avm\,,\\
h_2&\;=\; &h_1|ST^2\avp &-h_2|ST^{-2}\avm\,,
\end{alignat*}
imply $h_2=h_1|T^{-1}$. So they are equivalent to
\begin{equation}\label{eif1}
h \;=\; h|ST^3 \avp - h|T^{-1}ST^{-1}\avm \qquad \text{ on
}(-\gr^2,\gr)\,.
\end{equation} Applying $|(1-\nobreak T)$ to \eqref{eif1}, and using
\eqref{avTrel1}, we get
\[ h|(1-T) \;=\; h|(ST^3 - T^{-1}ST^{-1})\,,\]
as an identity in $\cV^\om_s(-\gr^2,\gr^{-1})$. Apply $|_{2s}T^{-1}$
to obtain as an equality in $\cV^\om_s(-\gr,\gr)$
\begin{equation}\label{4te|}
h|(1+ST^2) \;=\; h|(T^{-1}+T^{-1}ST^{-2})\,,
\end{equation}
which is~\eqref{4te}.
\end{proof}

%%%%%%%%%%%% Section 4 %%%%%%%%%%%%%%%

\section{Four term equation and cohomology} \label{D}
In the case of the three term equation \eqref{3te} the step from
solutions on $(0,\infty)$ to cohomology has been discussed in
\S\ref{sect-efMa}. In the case of the four term equation \eqref{4te},
the situation is more complicated, but the approach is essentially
based on the same ideas. The main steps are indicated in
\S\ref{sect-hypcoh}. In this section we prove Theorems~\ref{thm-th}
and~\ref{thm-par}.

\subsection{A model for parabolic and hyperbolic cohomology}
\label{D1}
Instead of the inhomogeneous cocycles used in \S\ref{B4} to describe
the first cohomology group, one can also employ homogeneous ones. To
an inhomogeneous cocycle $\psi$ corresponds the \emph{homogeneous
cocycle} $\tilde{c}: \Gmod^2\longrightarrow V$ given by
\begin{equation}
\label{D1.1}
\tilde c_{\gamma,\delta} = \psi_{\gamma\delta^{-1}}|\delta
\qquad(\g,\delta\in \Gmod)\,.
\end{equation}
satisfying for $\g,\delta,\e\in \Gmod$
\begin{equation}
\label{D1.2}
\tilde c_{\gamma\eps,\delta\eps} = \tilde c_{\gamma,\delta}|\eps
\quad \mbox{and} \quad \tilde c_{\gamma,\delta}+\tilde c_{\delta,\eps}
= \tilde c_{\gamma,\eps}.
\end{equation}
(This implies that $\tilde c_{\gamma,\gamma}=0$ and
$\tilde c_{\delta,\gamma}=-\tilde c_{\gamma,\delta}$.)
Coboundaries in $B^1(\Gmod;V)$ correspond to functions
\begin{equation}
\label{D1.3}
(\gamma,\delta)\mapsto v|\delta-v|\gamma
\end{equation}
with $v\in V$.

Homogeneous cocycles are parabolic if and only if by adding a
coboundary, we can arrange $\tilde c_{T,1}=0$.

A homogeneous parabolic cocycle satisfies
\begin{eqnarray}
\label{D1.4}
\tilde c_{T\gamma,\delta} &=& \tilde c_{T\gamma,\gamma}+\tilde
c_{\gamma,\delta} = 0|\gamma+\tilde c_{\gamma,\delta} \quad
\mbox{and} \\
\tilde c_{\gamma,T\delta} &=& \tilde c_{\gamma,\delta}+ \tilde
c_{1,T}|\delta = \tilde c_{\gamma,\delta}-0|\delta.
\end{eqnarray}
So $\tilde c$ induces a function
$(\Gamma_\infty\backslash \Gmod) ^2 \longrightarrow V$, where
$\Gamma_\infty$ is the subgroup of $\Gmod$ generated by $T$, which is
the subgroup fixing $\infty$. Conversely, every such function
satisfying \eqref{D1.2} induces a homogeneous parabolic cocycle.
Taking into account that the set $\Gamma_\infty\backslash\Gmod$ can
be identified with the set of cusps
$\PPP\Q =\Q\cup\infty \subset \PP$ via
$\gamma\mapsto \gamma^{-1}\infty$, we obtain the description of
$H^1_{\mathrm {par}} (\Gmod;V)$ as the quotient
$Z^1_{\PPP \Q}(\Gmod;V) / B^1_{\PPP\Q}(\Gmod;V)$, where
$Z^1_{\PPP \Q}(\Gmod;V)$ is the space of maps
$c:\PPP\Q\times\PPP\Q\longrightarrow V$ such that
$c_{\xi,\eta}+c_{\eta,\zeta} = c_{\xi,\zeta}$ and
$c_{\gamma^{-1}\xi,\gamma^{-1}\eta} = c_{\xi,\eta}|\gamma$ for
$\gamma\in \Gmod$ and $\xi,\eta,\zeta\in \PPP\Q$, and where
$B^1_{\PPP\Q}(\Gmod;V)$ consists of the subset of elements of the
form $c_{\xi,\eta} = f_\eta-f_\gamma$ with
$f:\PPP\Q\longrightarrow V$ satisfying
$f_{\gamma^{-1}\xi} = f_\xi|\gamma$.

\medskip

If we take another base point $\xi\in \PP$, we can work with cocycles
of the same type on other $\Gmod$-orbits in $\PP$. If the subgroup
$\Gmod_\x$ of $\Gmod$ leaving $\xi$ fixed is trivial, we get a
description of $H^1(\Gmod;V)$. In fact a homogeneous cocycle
$\tilde{c}$ on $\Gmod^2$ corresponds to a cocycle $c$ on
$\Gmod \xi\subset\PP$ by
\begin{equation}
\label{D1.6}
\tilde c_{\gamma,\delta} = c_{\gamma^{-1}\xi,\delta^{-1}\xi}.
\end{equation}

The situation is different if $\xi$ is a hyperbolic fixed point of $\Gmod$. Then the
procedure indicated above leads to the hyperbolic cohomology group
discussed in \S\ref{sect-hypcoh}. A corresponding homogeneous group
cocycle~$\ps: \g\mapsto \ps_\g =c_{\g^{-1}\x,\x}$ satisfies
\[ \psi_H \in V|(1-H)
\]
for a generator $H$ of the stabilizer $\Gmod_\xi$ of the point~$\xi$.
We will apply this for the stabilizer $\Gmod_{-\gr}$ for the golden
ratio $\gr=\frac{1+\sqrt 5}{2}$. A generator of $\Gmod_{-\gr}$ is
$TST^2$. Thus, we have a model of $H^1_\cH$. A group cocycle
corresponding $\ps$ to the cocycle $c$ on $\cH$ is given by
$\ps_\g=c_{\g^{-1}(-\gr),-\gr}$. In particular, $ c$ is determined by
$ c_{-\gr,\gr^{-1}} = - \ps_S$ and $ c_{-\gr,\gr} = - \ps_{T^{-1}S}$,
subject to the relations
\begin{equation}\label{hyprel}
 c_{-\gr,\gr^{-1}}|(1+S)\;=\;0\,,\qquad
 c_{-\gr,\gr}|(1+T^{-1}S+ST)\;=\;0\,;
\end{equation}
see \eqref{S-ST-rel}. Such a cocycle is a coboundary if there is
$v\in V$ such that $v|TST^2=v$ and
$c_{\g^{-1}(-\gr),\dt^{-1}(-\gr)} = v|\dt-v|\g$ for
$\g,\dt\in \Gmod$.

\subsection{Solutions of the four term equation and hyperbolic
cohomology} \label{D2}
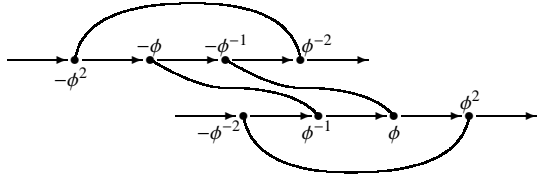
\begin{figure}
\[
\setlength{\unitlength}{1cm}
\begin{picture}(6,2.25)(-3,-2.25)
\put(-2.62,-.75){\circle*{.1}}
\put(-2.9,-1.1){$\scriptstyle -\gr^2$}
\put(-1.62,-.75){\circle*{.1}}
\put(-1.8,-.65){$\scriptstyle -\gr$}
\put(-.62,-.75){\circle*{.1}}
\put(-.9,-.65){$\scriptstyle -\gr^{-1}$}
\put(.38,-.75){\circle*{.1}}
\put(.4,-.65){$\scriptstyle \gr^{-2}$}
\put(-.38,-1.5){\circle*{.1}}
\put(-1,-1.8){$\scriptstyle -\gr^{-2}$}
\put(.62,-1.5){\circle*{.1}}
\put(.4,-1.8){$\scriptstyle \gr^{-1}$}
\put(1.62,-1.5){\circle*{.1}}
\put(1.55,-1.8){$\scriptstyle \gr$}
\put(2.62,-1.5){\circle*{.1}}
\put(2.5,-1.4){$\scriptstyle \gr^2$}
\put(-3.52,-.75){\vector(1,0){.8}}
\put(-2.52,-.75){\vector(1,0){.8}}
\put(-1.52,-.75){\vector(1,0){.8}}
\put(-.52,-.75){\vector(1,0){.8}}
\put(.48,-.75){\vector(1,0){.8}}
\put(2.72,-1.5){\vector(1,0){.8}}
\put(1.72,-1.5){\vector(1,0){.8}}
\put(.72,-1.5){\vector(1,0){.8}}
\put(-.28,-1.5){\vector(1,0){.8}}
\put(-1.28,-1.5){\vector(1,0){.8}}
\qbezier(-1.62,-.75)(-1,-1.125)(-.62,-1.125)
\qbezier(-.62,-1.125)(.3,-1.125)(.62,-1.5)
\qbezier(-.62,-.75)(0,-1.125)(.38,-1.125)
\qbezier(.38,-1.125)(1.3,-1.125)(1.62,-1.5)
\qbezier(-2.62,-.75)(-2.5,0)(-1.05,0)
\qbezier(-1.05,0)(.25,0)(.38,-.75)
\qbezier(2.62,-1.5)(2.5,-2.25)(1.05,-2.25)
\qbezier(1.05,-2.25)(-.25,-2.25)(-.38,-1.5)
\end{picture}
\]
\caption{Part of the $\Gmod$-orbit~$\cH$. Arrows denote the action
of~$T$, curves the action of~$S$. The horizontal coordinates of the
points correspond to their positions in~$\R$. Note that
$\gr^2=\gr+1$, $\gr^{-1}=\gr-1$ and $\gr^{-2}=2-\gr$.} \label{H-pict}
\end{figure}

In this subsection we prove Theorem~\ref{thm-th}. We work with
cocycles on $\cH=\Gmod(-\gr)$ as the model of hyperbolic cohomology.
For the computations we found the graph in
Figure~\ref{H-pict} useful.

\begin{lemma}\label{lem-coc-FE}Suppose that $c$ is a $1$-cocycle on
$\cH$ with values in $\cV^{\fs}_s$ such that
$\sing{c_{-\gr^{2},\gr}} \cap (-\gr^{2},\gr)=\emptyset$. Then the
restriction of $c_{-\gr^{2},\gr}$ to $(-\gr^{2},\gr)$ is in
$\FEv_s(-\gr^{2},\gr)_\om$.
\end{lemma}
\begin{proof}Cocycles on $\cH$ satisfy the relations in~\eqref{D1.2}.
Hence
\begin{align}\label{4tev}
c_{-\gr^{2},\gr}|(1-T^{-1}ST^{-2}) &\;=\;
c_{-\gr^{2},\gr}-c_{\gr^2,\gr} \;=\; c_{-\gr^2,\gr^2} \\
\nonumber
& \;=\; c_{-\gr^{2},-\gr} + c_{-\gr,\gr^2} \;=\;
c_{-\gr^2,\gr}|(T^{-1}-ST^2)\,.
\end{align}
The fact that the restriction $g$ of $c_{-\gr^2,\gr}$ to
$(-\gr^2,\gr)$ has no singularities in $(-\gr^2,\gr)$ implies that
relation \eqref{4tev} holds on $(-\gr,\gr)$ as an identity of real
analytic functions. Hence $g\in \FEv_s(-\gr^2,\gr)_\om$.
\end{proof}

\begin{proposition}\label{prop-thh}Each element
$g\in \FEv_s(-\gr^2,\gr)_\om$ is the restriction
\[ c_{-\gr^2,\gr}|_{(-\gr^2,\gr)} \]
for a unique $\Gmod$-cocycle $c$ on $\cH$ with values in
$\cV^{\fs}_s$. This cocycle satisfies
\[\sing{ c_{-\gr^2,\gr} } \subset \{-\gr^2,\gr\}\,,\;\;\;
\sing{c_{-\gr,\gr^{-1}}}\subset \{-\gr,\gr^{-1}\}\,,\;\;\;
\sing{c_{-\gr,\gr} }\subset\{-\gr,\gr\}\,. \]
The induced map $\th:\FEv_s(-\gr^2,\gr)_\om) \longrightarrow
H^1_\cH(\Gmod;\cV^{\fs}_s)$ is injective.
\end{proposition}
\begin{proof}First we assume that a cocycle $c$ exists such that $g$
is equal to the restriction of $c_{-\gr^2,\gr}$ to $(-\gr^2,\gr)$.
The cocycle relations imply
\begin{alignat*}2 g|ST^2 &\;=\; c_{-\gr,-\gr^2}&\quad\text{ on
}&(-\gr,-\gr^2)_c\,,
\displaybreak[0]\\
g|(1+ST^2)&\;=\; c_{-\gr^2,\gr}+c_{-\gr,-\gr^2}\;=\;
c_{-\gr,\gr}&\text{ on }&(-\gr,\gr)\,,
\displaybreak[0]\\
g|(1+ST^2)STS&\;=\; c_{\gr^2,\gr^{-1}}&\text{ on
}&(\gr^2,\gr^{-1})_c\,,
\displaybreak[0]\\
g|T^{-1}&\;=\; c_{-\gr,\gr^2}&\text{ on }&(-\gr,\gr^2)\,,
\displaybreak[0]\\
g|\left( T^{-1}+(1+ST^2)STS\right)&\;=\;
c_{-\gr,\gr^2}+c_{\gr^2,\gr^{-1}} \;=\; c_{-\gr,\gr^{-1}}&\text{ on
}&(-\gr,\gr^{-1})\,.
\end{alignat*}
This shows that the restriction of $c_{-\gr,\gr^{-1}}$ to
$(-\gr,\gr^{-1})$ is determined by $g$. {}From \eqref{hyprel} we know
that $c_{-\gr,\gr^{-1}}|S = - c_{-\gr,\gr^{-1}}$. Hence $g$
determines the restriction of $c_{-\gr,\gr^{-1}}$ to
$(\gr^{-1},-\gr)_c$ as well. So $g$ determines $c_{-\gr,\gr^{-1}}$ as
an element of~$\cV^{\fs}_s$.

The situation for the other generator $c_{-\gr,\gr}$ is slightly more
complicated.
\begin{alignat*}2 c_{-\gr,\gr} &\;=\; g|(1+ST^2) &\quad\text{ on
}&(-\gr,\gr)\,,\displaybreak[0]\\
c_{\gr,\gr^{-2}}&\;=\; c_{-\gr,\gr}|T^{-1}S \;=\;
g|(1+ST^2)|T^{-1}S&\quad\text{ on }&
(\gr,\gr^{-2})_c\,,\displaybreak[0]\\
c_{\gr^{-2},-\gr}&\;=\; c_{-\gr,\gr}|ST\;=\; g|(1+ST^2)|ST&\quad\text{
on }&(\gr^{-2},-\gr)_c\,.
\end{alignat*}
The relation $c_{-\gr,\gr}+c_{\gr,\gr^{-2}}+c_{\gr^{-2},-\gr}=0$
implies that $c_{-\gr,\gr}$ is determined by $g$ on each of the
cyclic intervals $(-\gr,\gr^{-2})$, $(\gr^{-2},\gr)$ and
$(\gr,-\gr)_c$. So any $\cV^\fs_s$-valued cocycle $c$ on $\cH$ is
determined by the restriction of $c_{-\gr^2,\gr}$ to $(-\gr^2,\gr)$.

This reasoning also shows how to
construct $c$ from~$g$. We put for a given function
$g\in \FEv_s(-\gr^2,\gr)_\om$:
\begin{alignat}2 \label{hkdef}
h&\;=\;& g|(1+ST^2)&\in \cV^\om_s(-\gr,\gr)\,,\\
\nonumber
k&\;=\;& g|T^{-1}+h|STS&\in \cV^\om_s(-\gr,\gr^{-1})\,.
\end{alignat}
By the reasoning given above we should have:
\begin{align}\label{ccdef}
c_{-\gr,\gr^{-1}} &\;:=\; \begin{cases}
k&\text{ on }(-\gr,\gr^{-1})\,,\\
-k|S&\text{ on }(\gr^{-1},-\gr)_c\,,
\end{cases}\\
\nonumber
c_{-\gr,\gr}&\;:=\; \begin{cases}
h&\text{ on }(-\gr,\gr)\,,\\
-h|(T^{-1}S+ST)&\text{ on }(\gr,-\gr)_c\,.
\end{cases}
\end{align}
To see that this
indeed defines a cocycle, we choose the base point $-\gr$ and
consider not the potential cocycle $c$ on $\cH$, but the
corresponding cocycle $\ps$ on~$\Gmod$:
\[ \ps_S \;=\; -c_{-\gr,\gr^{-1}}\,,\qquad \ps_{T^{-1}S} \;=\; -
c_{-\gr,\gr}\,.\]
The relations \eqref{S-ST-rel} turn out to be satisfied. So indeed
there exists a cocycle $c$ as desired.

For the singularities, we note first that the expressions above for
$c_{-\gr,\gr}$ and $c_{-\gr,\gr^{-1}}$ in terms of~$g$ imply that
\[ \sing{c_{-\gr,\gr^{-1}}} \subset \{-\gr,\gr^{-1}\}\,,\qquad \sing{
c_{-\gr,\gr}} \subset\{-\gr,\gr\}\,.\]

The cocycle $c$ is determined by $c_{-\gr,\gr^{-1}}$ and
$c_{-\gr,\gr}$. We now check that the restriction of $c_{-\gr^2,\gr}$
to $(-\gr^2,\gr)$ is equal to~$g$, as desired, and cannot have
singularities in $(\gr,-\gr^2)_c$. The cocycle relations imply
that
\begin{equation}
\label{cr} c_{-\gr^2,\gr}=c_{-\gr,\gr^{-1}}|T+
c_{-\gr,\gr}+c_{-\gr,\gr}|ST \,.
\end{equation}
On various intervals we have:
\[ \renewcommand\arraystretch{1.5}
\begin{array}{|c|c|c|c|c|}\hline
&(-\gr^2,-\gr)&(-\gr,-\gr^{-2})&(-\gr^{-2},\gr)&(\gr,-\gr^2)_c\\
\hline
c_{-\gr,\gr^{-1}}|T& k|T&k|T&-k|ST&-k|ST \\
c_{-\gr,\gr} &-h|(T^{-1}S+ST)& h&h&-h|(T^{-1}+ST)\\
c_{-\gr,\gr}|ST& h|ST&-h|(1+T^{-1}S)&h|ST&h|ST \\ \hline
c_{-\gr^2,\gr}& g & g & g
& -k|ST-h|T^{-1}
\\ \hline
\end{array}
\]
On $(-\gr^{-2},\gr)$ we use the four term equation:
\begin{align*}
c_{-\gr^2,\gr}&\;=\; - g|T^{-1}ST + h|(-ST^2 + 1+ST)\\
&\;=\; g | \left( -T^{-1}ST + (1+ST^2)(1+ST-ST^2) \right)\\
&\!\!\stackrel{\eqref{4te|}}=\; - g|(1+ST^2-T^{-1}ST^{-2})ST +
g|\left( 1-ST^2ST^2+ST+ST^2ST\right)\\
&\;=\; g|\left( T^{-1}ST^{-2}ST+1-ST^2ST^2\right)
\;\stackrel{\eqref{B1.1}}=\; g\,.
\end{align*}
Hence $c_{-\gr^2,\gr}= g$ on $(-\gr^2,-\gr) \cup (-\gr,-\gr^{-2})
\cup( -\gr^{-2},\gr)$. Since $g$ is analytic on~$(-\gr^2,\gr)$, the
points $-\gr$ and $-\gr^{-2}$ are not singularities of
$c_{-\gr^2,\gr}$, by the definition of $\cV^{\fs}_s$ as an inductive
limit. Furthermore $c_{-\gr^2,\gr}$ is given by the analytic function
$-k|ST-h|T^{-1}$ on $(\gr,-\gr^2)_c$. This shows that
$\sing {c_{-\gr^2,\gr}} \subset\{-\gr^2,\gr\}$.

To show that the map
$\th: \FEv_s(-\gr^2,\gr)_\om \longrightarrow H^1_\cH(\Gmod;\cV^\fs_s)$
given by $\th: g \mapsto [c]$ is injective, we check that the cocycle
$c$ corresponding to $g$ can be a coboundary only if $g=0$. If $c$ is
a coboundary, then $c_{\g^{-1}(-\gr),\dt^{-1}(-\gr)} = v|\dt - v|\g$.
In particular, $v|TST^2$ should be equal to~$v$ for some $v\in
\cV^\fs_s$.

We use that $TST^2$ is a hyperbolic element of~$\Gmod$ fixing $-\gr$
and $\gr^{-1}$. Conjugation in $\PSL\R$ transforms it to
$\eta=\Matrix {\gr^2}00{\gr^{-2}}$, fixing $0$ and $\infty$.
Let $w\in \cV^{\fs}_s$ be invariant
under~$\eta$. The action of $\eta$ on $(0,\infty)$ and $(-\infty,0)$
is by $x\mapsto  \gr^4 x$. So if $w$ were to have singularities in
$\R\setminus\{0\}$, then there would be infinitely many,
contradicting the definition of $\cV^{\fs}_s$.
On $\R$ we have
$w|\eta(x) = \gr^{4s} w(\gr^4x)$. Inserting this into a power series
expansion converging on a neighborhood of~$0$, we
see that if $w$ is analytic
at~$0$ it vanishes. The same holds at~$\infty$. So if $w\neq 0$, then
$\sing w = \{0,\infty\}$.

Conjugating back, we see that if $c$ is a non-zero coboundary, then it
has the form $c_{\g^{-1}(-\gr),\dt^{-1}(-\gr)} = v|(\dt-\nobreak\g)$,
where $v\in \cV^{\fs}_s$ satisfies $v|TST^2=v$ and
$\sing v = \{-\gr,\gr^{-1}\}$. We consider the singularities of
$c_{-\gr^2,\gr} = c_{T^{-1}(-\gr),TS(-\gr)} = v|ST^{-1}-v|T$. Now
$\sing{v|ST^{-1}}=\{\gr,-\gr^{-1}\}$, and
$\sing{v|T} = \{-\gr^2, -\gr^{-2}\}$. So $c_{-\gr^2,\gr}$ has
singularities at $-\gr^{-1}$ and $-\gr^{-2}$, in contradiction to the
analyticity of~$g$ on $(-\gr^2,\gr)$.
\end{proof}

We note that $H^1_\cH(\Gmod;\cV^{\fs}_s)$ is a subspace of
$H^1(\Gmod;\cV^{\fs}_s)$. The inclusion
$\cV^\om_s\hookrightarrow\cV^{\fs}_s$ induces a natural map
$H^1(\Gmod;\cV^\om_s) \longrightarrow  H^1(\Gmod;\cV^{\fs}_s)$.

\begin{proposition}\label{prop-om}The subspace
$\th(\FEv_s(-\gr^2,\gr)_\om) $ of $ H^1(\Gmod;\cV^{\fs}_s)$ is
contained in the image of $H^1(\Gmod;\cV^\om_s)$ in
$ H^1(\Gmod;\cV^{\fs}_s)$.
\end{proposition}
\begin{proof}Let $g\in \FEv_s(-\gr^2,\gr)_\om$, and let $c$ be the
cocycle on~$\cH$ representing $\th(\g)$. Since
$\sing{c_{-\gr,\gr^{-1}}} \subset \{ -\gr,\gr^{-1}\}$,
Proposition~\ref{D2.6} implies that there are $K_{-\gr}$ and
$K_{\gr^{-1}}$ in $\cV^{\fs}$ with $\sing{K_{-\gr}}\subset\{-\gr\}$
and $\sing{K_{\gr^{-1}}}\subset\{\gr^{-1}\}$ such that
$c_{-\gr,\gr^{-1}} = K_{-\gr}-K_{\gr^{-1}}$. Since
$c_{-\gr,\gr^{-1}}|(1+\nobreak S)=0$, we have
$K_{\gr^{-1}} = K_{-\gr}|S + \al$ for some $\al\in \cV^\om_s$.
Similarly, there are $H_{-\gr}, H_\gr\in \cV^{\fs}_s$ with
$\sing{H_{-\gr}}\subset\{-\gr\}$ and $\sing{H_\gr}\subset \{\gr\}$
such that $c_{-\gr,\gr} = H_{-\gr}-H_\gr$. There exists a function
$\bt \in \cV^\om_s$ such that $H_\gr = H_{-\gr}|T^{-1}S+\bt$.

We have, with \eqref{cr}:
\[c_{-\gr^2,\gr} \;\in\;c_{-\gr,\gr^{-1}}|T + c_{-\gr,\gr}|(1+\nobreak
ST)
\in K_{-\gr}|(T-\nobreak ST) + H_{-\gr}|(ST-\nobreak T^{-1}S) +
\cV^\om_s \,.\]
Since $\sing{c_{-\gr^2,\gr}}\subset\{-\gr^2,\gr\}$, we conclude from
the singularities of the various terms that
$K_{-\gr}|ST \in H_{-\gr}|ST +\cV^\om_s$. Hence
\begin{equation} K_{-\gr}-H_{-\gr}\in \cV^\om_s\,.\end{equation}
\smallskip

The class $[c]\in H^1_\cH(\Gmod;\cV^{\fs}_s)$ considered as a class in
$H^1(\Gmod;\cV^{\fs}_s)$ is given by the group cocycle
\[ \ps_S \;=\; -c_{-\gr,\gr^{-1}}\,,\qquad \ps_{T^{-1}S} \;=\; -
c_{-\gr,\gr}\,.\]
We add to it the coboundary $\g\mapsto K_{-\gr}|(1-\nobreak\g)$,
obtaining a cocycle $\tilde\ps$ in the same class. It satisfies
\begin{align}
\label{tpsi-on-gen}
\tilde\ps_S &\;=\;
-K_{-\gr}+K_{\gr^{-1}}+K_{-\gr}|(1-S) \;=\; \al\;\in\; \cV^\om_s\,,\\
\nonumber
\tilde\ps_{T^{-1}S} &\;=\;
-H_{-\gr}+H_\gr + K_{-\gr}|(1-\nobreak T^{-1}S) \\
\nonumber
&\;=\; (K_{-\gr}-\nobreak H_{-\gr})|(1-\nobreak T^{-1}S) +\bt \;\in\;
\cV^\om_s\,.
\end{align}
Thus the class $[\ps] = [\tilde\ps]$ in $H^1(\Gmod;\cV^\fs_s)$ is the
image of the class
$[\tilde\ps]\in H^1(\Gmod;\cV^\om_s)$.
\end{proof}

\begin{proposition}\label{prop-thsurj}The subspace
$\th\left( \FEv_s(-\gr^2,\gr)_\om \right)$ of
$H^1(\Gmod;\cV^{\fs}_s)$ is equal to the image of
$H^1(\Gmod;\cV^\om_s)$.
\end{proposition}
\begin{proof} We use a hyperbolic one-sided average, as discussed at
the end of~\S\ref{C1}. Starting with a cocycle
$\tilde \ch\in Z^1(\Gmod;\cV^\om_s)$ we define
\[ A= \tilde \ch_{TST^2} | \av_{TST^2}^+ \in \cV^{\fs}_s\,,\]
with $\sing A \subset \{-\gr\}$ and
$A|(1-\nobreak TST^2) = \tilde \ch_{TST^2}$. We have used that $-\gr$
is the repelling fixed point of $TST^2$.

Now $\ch:\g\mapsto \tilde\ch_\g - A|(1-\nobreak\g)$ defines a
$\cV^{\fs}_s$-valued cocycle in the same cohomology class in
$H^1(\Gmod;\cV^{\fs}_s)$ as $\tilde \ch$. It satisfies
$\sing{\ch_\g}\subset \{-\gr,\g^{-1}(-\gr)\}$. Moreover,
$\ch_{TST^2}=0$, so $\ch$ is hyperbolic for the conjugacy class of
$TST^2$. Hence it corresponds to a cocycle $c$ on~$\cH$, such that
$c_{-\gr^2,\gr} = A|T - A|ST^{-1}-\tilde\ch_{ST^{-2}}|T$ has
singularities at most in $\{-\gr^2,\gr\}$. This means that the class
of $\tilde \ch$ is in the image of~$\th$.

Starting from this cocycle $c$, one can take $K_{-\gr}=H_{-\gr} = A$,
and get back $\tilde\ps = \tilde\ch$.
\end{proof}
The results in this subsection imply Theorem~\ref{thm-th}.

\subsection{Eigenfunctions of the transfer operator and
cohomology}\label{sect-eif-coh}In this final section we prove
Theorem~\ref{thm-par}.

We use the same notations as in the previous subsection: $g\in
\FEv_s(-\gr^2,\gr)_\om$, $c\in Z^1_\cH(\Gmod;\cV^{\fs}_s)$,
$K_{-\gr}\in \cV^\fs_s$ with $\sing{K_{-\gr}} \subset\{-\gr\}$,
$\ps\in Z^1(\Gmod;\cV^{\fs}_s)$ and
$\tilde\ps\in Z^1(\Gmod;\cV^\om_s)$, related by:
\begin{align}
\label{g-c}
g&\;=\; c_{-\gr^2,\gr}\text{ restricted to }(-\gr^2,\gr)\,,\\
\label{c-ps}
c_{\g^{-1}(-\gr),\dt^{-1}(-\gr)} &\;=\; \ps_{\g\dt^{-1}}|\dt,\qquad
\ps_\g \;=\;
-c_{-\gr,\g^{-1}(-\gr)}\,,\\
\label{ps-tps}
\tilde\ps_\g&\;=\; \ps_\g + K_{-\gr}|(1-\g)\,,\\
\label{K-ps}
K_{-\gr}&\;=\; \tilde\ps_{TST^2}|\av_{TST^2}^+\,.
\end{align}

We consider $P\in \cV^\om_s(-\gr^2,\gr)$ given by
\begin{equation}
P\;=\; g|ST^3\avp-g|T^{-1}ST^{-1}\avm-g\,.
\end{equation}
 So $P$ measures how much $(g,g|T^{-1})$ differs from an
$1$-eigenfunction of $\LL$ in~\eqref{LL-f}. Application of
\eqref{avTrel1} shows that $P|T=P$ on $(-\gr^2,\gr^{-1})$. Hence the
periodic function $P$ extends as an element of
$\cV^\om_s(\R)^T \subset\cV^{\fs}_s$. To obtain a representation of
$P$ on $\R$, we use \eqref{g-c}--\eqref{ps-tps} to obtain on
$(-\gr^2,\gr)$:
\begin{align*}
P&\;=\; \ps_{TST}|T^{-1}S \left( ST^3\avp-T^{-1}ST^{-1}\avm-1\right)\\
&\;=\; \tilde\ps_{TST}|\left( T^2\avp-S\avm -T^{-1}S\right)\\
&\qquad\hbox{} \nobreak
+ K_{-\gr}|\left( \left(TST^3 - T^2\right)\avp + \left(
S-ST^{-1}\right)
\avm + T^{-1}S-T\right)\,.
\end{align*}
{}From \eqref{K-ps} it follows that 
$K_{-\gr} |TST^2 = K_{-\gr} - \tilde\ps_{TST^2}$. With
\eqref{avTrel2} we find for the contribution of~$K_{-\gr}$:
\begin{align*}
K_{-\gr}|&(T-T^2)\avp - \tilde\ps_{TST^2}|T\avp
-K_{-\gr}|ST^{-1}(1-T)\avm
+ K_{-\gr}|(T^{-1} S-T)\\
&\;=\; K_{-\gr}|\left( T - ST^{-1} +T^{-1}S-T\right)
-\tilde\ps_{TST^2}|T\avp\\
&\;=\; \left( K_{-\gr}|TST^2 + \tilde\ps_{TST^2}\right)|T^{-1}S -
K_{-\gr}|ST^{-1} - \tilde\ps_{TST^2}|T\avp\\
&\;=\; \tilde\ps_{TST^2}|(T^{-1}S-T\avp)\,.
\end{align*}
This leads to an expression valid on~$\R$:
\begin{align}
\nonumber
P&\;=\; \tilde\ps_{TST}|\left( T^2 \avp-S\avm -T^{-1}S +
T(T^{-1}S-T\avp)\right)
\\
\nonumber
&\qquad\hbox{}
+ \tilde\ps_T| (T^{-1}S +1 -\avp)\\
\label{P-glob}
&\;=\; \tilde\ps_{TST}|(S-T^{-1}S) + \tilde\ps_T|(T^{-1}S+1)
- \tilde\ps_T |\avp -\tilde\ps_{TST}|S\avm\,.
\end{align}

\begin{proof}[Proof of Theorem~\ref{thm-par}] Suppose that $P=0$,
i.e., $g\in \FEv_s(-\gr^2,\gr)_\om$ corresponds to an eigenfunction
$(g,g|T^{-1})$ of $\LL$. We note that
$\tilde \ps_{TST}, \tilde\ps_T \in \cV^\om_s$, and apply
 Lemma~\ref{C1.8} to conclude that $\tilde\ps_T
|\avp = \tilde\ps_T|\avm \in \cV^{\smp}_s$, hence
$\tilde \ps_T \in \cV^{\smp}_s|(1-\nobreak T)$, which is the
condition for a cocycle to be parabolic. Thus we conclude that
$[\tilde\ps] \in H^1_\parb(\Gmod;\cV^\om_s,\cV^{\smp}_s)$.

Conversely, suppose that
$\tilde\ps\in Z^1_\parb(\Gmod;\cV^\om_s,\cV^{\smp}_s)$.
 By parabolicity $\tilde\ps_T \in \cV^{\smp}_s|(1-\nobreak T)$, and
hence $\tilde\ps_T|\avp=\tilde\ps_T|\avm$ by Lemma~\ref{8.17}. In
\eqref{P-glob} we see that $P$ is in $\cV^\om_s + \cV^\om_s|\avm$.
Hence $P(x)$ has an asymptotic behavior as $x\downarrow-\infty$ of
the form indicated in \eqref{avm-as}. Since $P$ is also periodic, it
vanishes.
\end{proof}

%%%%%% References 

\bibliographystyle{amsalpha}

\end{document}